\documentclass[a4paper,10pt,leqno,twoside]{article}
\usepackage{expl3}
\usepackage{amssymb,amsmath,amsthm}
\usepackage{stix}
\DeclareSymbolFont{calletters}{OMS}{cmsy}{m}{n}
\DeclareSymbolFontAlphabet{\mathcal}{calletters}
\usepackage{tgpagella}
\usepackage[T1]{fontenc}
\usepackage[utf8]{inputenc}
\usepackage[english]{babel}
\usepackage{authblk}
\usepackage{enumitem}
\usepackage{stmaryrd}
\usepackage{verbatim}

\usepackage[
backend=bibtex,
maxnames=50,
eprint=false,
sorting=nyt,
giveninits=true,
isbn=false,
url=false
]{biblatex}
\DeclareFieldFormat*{title}{\mkbibemph{#1}}
\DeclareFieldFormat*{citetitle}{\mkbibemph{#1}}
\DeclareFieldFormat{journaltitle}{#1}
\renewbibmacro*{in:}{
	\ifentrytype{article}
	{}
	{\printtext{\bibstring{in}\intitlepunct}}}
\newbibmacro*{pubinstorg+location+date}[1]{
	\printlist{#1}
	\newunit
	\printlist{location}
	\newunit
	\usebibmacro{date}
	\newunit}
\renewbibmacro*{publisher+location+date}{\usebibmacro{pubinstorg+location+date}{publisher}}
\renewbibmacro*{institution+location+date}{\usebibmacro{pubinstorg+location+date}{institution}}
\renewbibmacro*{organization+location+date}{\usebibmacro{pubinstorg+location+date}{organization}}

\makeatletter
\newcommand\blfootnote[1]{\xdef\@thefnmark{}\@footnotetext{\hspace*{-1.8em}#1}}
\makeatother

\usepackage{tikz,tikz-cd}
\tikzcdset{arrow style=math font}
\newlength{\perspective}
	\title{Zariski-constructible functions of class $\mathcal{C}^k$}
	\date{November 23, 2026}

	\author{Jean-Baptiste Campesato\thanks{jb.campesato@univ-angers.fr}}
	\author{Jean-Philippe Monnier\thanks{jean-philippe.monnier@univ-angers.fr}}
	\author{Enrico Savi\thanks{enrico.savi@univ-angers.fr}}
	\affil{\normalsize\itshape Univ Angers, CNRS, LAREMA, SFR MATHSTIC, F-49000 Angers, France}
	
    \usepackage[left=3.5cm,right=3.5cm,top=4cm,bottom=4cm,includeheadfoot,headsep=.5cm]{geometry}
	\usepackage{fancyhdr}
	\renewcommand{\thepage}{\arabic{page}}
	\usepackage{hyperref}
	\usepackage{color}
	\definecolor{darkred}{rgb}{.5,0,0}
	\definecolor{darkgreen}{rgb}{0,.5,0}
	\definecolor{darkblue}{rgb}{0,0,.5}
	\hypersetup{
		pdfencoding=unicode,
		colorlinks=true,
		linkcolor=darkred,
		citecolor=darkgreen,
		urlcolor=darkblue,
		bookmarks=false,
		pdftoolbar=true,
		pdfmenubar=true,
		pdffitwindow=true,
		pdfstartview={FitH},
		pdfcreator={},
		pdfproducer={},
		pdfnewwindow=true,
		pdftitle={Zariski-constructible functions of class C^k},
		pdfauthor={Jean-Baptiste Campesato, Jean-Philippe Monnier, Enrico Savi},
	}
	
	\newtheorem*{thm*}{Theorem}
	\newtheorem{thm}{Theorem}[section]
	\newtheorem{lemma}[thm]{Lemma}
	\newtheorem{prop}[thm]{Proposition}
	\newtheorem{cor}[thm]{Corollary}
	\newtheorem*{question}{Question}
	\theoremstyle{definition}
	\newtheorem{defn}[thm]{Definition}
	
	\newtheorem{notation}[thm]{Notation}
	\newtheorem{eg}[thm]{Example}
	\theoremstyle{remark}
	\newtheorem{rem}[thm]{Remark}
	\newtheorem*{claim}{Claim}
	\renewcommand\qedsymbol{\ensuremath{\blacksquare}}
	\newcommand\cqed{{\renewcommand\qedsymbol{\ensuremath{\square}}\qed}}
	
	\newcommand{\ZC}{\mathcal{ZC}}
	\newcommand{\ZCc}{\mathcal{ZC}^0}

	\newcommand{\lk}{\operatorname{Lk}}
	\newcommand{\st}{\operatorname{St}}
	\newcommand{\dst}{\partial\operatorname{St}}
	\newcommand{\relint}{\operatorname{relint}}
	\newcommand{\zcl}[1]{\overline{#1}^{\ZC}}
	\newcommand{\zarcl}[1]{\overline{#1}^{\text{Zar}}}
	\newcommand{\zbd}{\partial_{\ZC}}
	\newcommand{\ebd}{\partial_{\text{Euc}}}
	\newcommand{\ecl}[1]{\overline{#1}^{\text{Euc}}}
\newcommand{\eint}[1]{\mathrm{Int}_{\text{Euc}}(#1)}
\newcommand{\zint}[1]{\mathrm{Int}_{\ZC}(#1)}
	
	\newcommand{\I}{\operatorname{\mathcal I}}
	\newcommand{\colim}{\underrightarrow{\lim}}

\begin{document}
		\maketitle
		
		\blfootnote{\ \\During the preparation of this work, the authors were supported by the ANR project NewMIRAGE (ANR-23-CE40-0002). The third author was partially supported by GNSAGA of INdAM. \\
			
			
		}
		
		\begin{abstract}
			Let $R$ be a real closed field. A Zariski-constructible subset of $R^n$, $\mathcal{ZC}$-set for short, is a Boolean combination of algebraic subsets of $R^n$. Euclidean closed Zariski-constructible subsets of $R^n$ are the closed subsets of a Noetherian topology finer than the Zariski topology of $R^n$.
		
			Let $k\in\mathbb{N}\cup\{\infty\}$. Let $U\subset R^n$ be a $\mathcal{ZC}$-open set. We say that a function $f:U\to R$ of class $\mathcal{C}^k$ is a Zariski-constructible function, $\mathcal{ZC}^k$-function for short, if its graph is a $\mathcal{ZC}$-set. For instance, regular and $k$-regulous functions on $R^n$ are $\mathcal{ZC}^k$-functions for any $k\in\mathbb{N}\cup\{\infty\}$. We show that the ring $\mathcal{ZC}^k(U)$ of $\mathcal{ZC}^k$-functions on a $\mathcal{ZC}$-open set is way richer, even for $k=\infty$, although $\infty$-regulous functions are just regular. Then, we extend the notion of $\mathcal{ZC}^k$-function over any $\mathcal{ZC}$-locally closed subset $X$ of $R^n$.
			
			This paper is mainly devoted to the study of the geometric and algebraic properties of the ring $\mathcal{ZC}^k(X)$ of $\mathcal{ZC}^k$-functions on a $\mathcal{ZC}$-locally closed subset $X$ of $R^n$. From now on, assume that $k\in\mathbb N$. We show that $\mathcal{ZC}$-closed subsets $X$ of $R^n$ are exactly the zero loci of $\mathcal{ZC}^k$-functions on $R^n$. We show that the Nullstellensatz holds in the ring $\mathcal{ZC}^k(U)$, where $U$ is a $\mathcal{ZC}$-open subset of $R^n$, and that this ring is radically principal. As a consequence of the Nullstellensatz, we prove that any $\mathcal{ZC}^k$-function on a $\mathcal{ZC}$-locally closed subset $X$ of $R^n$ extends to a $\mathcal{ZC}$-open neighborhood $U$ of $X$. Then, we translate our results into the language of sheaves of $\mathcal{ZC}^k$-functions on $\mathcal{ZC}$-locally closed sets, from which we deduce a $\mathcal{ZC}^k$-version of Cartan's theorems A and B.
		\end{abstract}
		
		\newpage
		
		\tableofcontents
		
		\section{Introduction}
		
		\subsection{Context}
		
	Real algebraic geometry studies, among other things, real algebraic sets, that is, the solution sets of polynomial equations over $\mathbb{R}$, or more generally, over a real closed field $R$. There are two different approaches in this area which have emerged in order to pursue different questions in the field.\\
	
	The first approach consists of including real algebraic geometry as a special case of complex algebraic geometry. More precisely, let $X(\mathbb C)\coloneq\{x\in \mathbb{C}^n\,:\,p_1(x)=\dots=p_k(x)=0\}$ be the complex algebraic subset of $\mathbb{C}^n$ defined over $\mathbb R$ by the polynomials $p_1,\dots,p_k\in\mathbb{R}[x_1,\dots,x_n]$. In this case, an equivariant version of Hilbert's Nullstellensatz gives a correspondence between algebraic subsets of $\mathbb C^n$ defined over $\mathbb{R}$ and radical ideals of $\mathbb{R}[x_1,\dots,x_n]$. Equivalently, in modern terms, this gives a correspondence between algebraic subsets of $\mathbb C^n$ defined over $\mathbb R$ and Zariski closed subsets of $\operatorname{Spec}(\mathbb R[x_1,\dots,x_n])$ and therefore modern techniques of algebraic geometry relying on the theories of sheaves and schemes find natural application. This approach has proven to be very fruitful for the study of complex algebraic sets defined over $\mathbb{R}$ from a birational point of view but also when looking for constraints on the topology of the real algebraic set $X(\mathbb R)\coloneq\{x\in \mathbb{R}^n\,:\,p_1(x)=\dots=p_k(x)=0\}$ in terms of the topology of the complex algebraic set $X(\mathbb{C})\subset \mathbb{C}^n$. However, when trying to characterize all possible topological types of real algebraic sets $X(\mathbb R)\subset \mathbb{R}^n$, both for singular and nonsingular ones, a different perspective is often preferable, we refer for instance to the monographs \cite{AK92,BCR}.\\
	
	The second perspective consists of developing real intrinsic techniques to study the geometry and the topology of the real algebraic set $X(\mathbb R)\subset \mathbb{R}^n$, which is the one we follow in the present paper. As the main objects of study are real algebraic sets, we will simply say that $X\subset \mathbb R^n$ is an algebraic set if there are some polynomials $p_1,\dots,p_k\in\mathbb{R}[x_1,\dots,x_n]$ such that $X=\{x\in\mathbb R^n\,:\,p_1(x)=\dots=p_k(x)=0\}$. Observe that every real algebraic set $X$ as above comes with an associated complex algebraic set usually called the complexification $X_{\mathbb C}$ of $X$ defined as the Zariski closure of $X$ in $\mathbb{C}^n$. Observe that $X_{\mathbb C}\subset \mathbb C^n$ is an algebraic set defined over $\mathbb{R}$ but it does not always coincide with the complex algebraic set $\{x\in\mathbb{C}^n\,:\,p_1(x)=\dots=p_k(x)=0\}$ if $X=\{x\in\mathbb{R}^n\,:\,p_1(x)=\dots=p_k(x)=0\}$. As an example, observe that $\varnothing=\{x\in \mathbb R\,:\,x^2+1=0\}$ and its complexification is $\varnothing_{\mathbb{C}}=\varnothing\neq\{\pm i\}=\{x\in \mathbb C\,:\,x^2+1=0\}$. The complexification $X_\mathbb{C}$ of an algebraic set $X$ is a very useful object in order to apply techniques of complex algebraic geometry to deduce properties of $X$. However, the previous simple example shows that the Nullstellensatz fails, as the ideal $(x^2+1)$ of $\mathbb R[x]$ is maximal, its real solution set is empty but $(x^2+1)\neq \mathbb R[x]$. Nevertheless, the failure of the Nullstellensatz shows as well that the ring of regular functions $\mathcal{R}(X)$ on $X\subset\mathbb R^n$ is particularly rich; as an example, the function $f:\mathbb{R}\to \mathbb{R}$ defined as $f(x)\coloneq\frac{1}{x^2+1}$ is regular over the whole $\mathbb R$ despite not being regular over the complexification $\mathbb C$ of $\mathbb R$, as it has two complex poles. The use of such functions has played a fundamental role in the study of the topological types of real algebraic sets after the pioneering work of Nash \cite{Nash52} and the later discovery of Tognoli \cite{Tog73} that every compact $\mathcal{C}^\infty$ manifold is diffeomorphic to a nonsingular real algebraic set, verifying a conjecture posed by Nash. We refer the interested reader again to the monographs \cite{AK92,BCR} for additional information.\\
	
	Although regular functions on real algebraic sets play such an important role in real algebraic geometry, they show some limitations from an algebraic point of view; the most important ones are again the failure of the Nullstellensatz in this ring and the failure of Cartan's Theorem A for quasi-coherent sheaves of modules on $X$ endowed with the structure sheaf of regular functions $\mathcal{R}_X$ on $X$, see \cite[Ex.\,12.1.5]{BCR}. Several classes of real functions, containing regular ones, have then been introduced in order to overcome these issues. All the classes considered below consist of semialgebraic functions.\\

		In 1988, Kurdyka introduced arc-symmetric sets and arc-analytic functions in \cite{Kur88}. A semialgebraic set $X\subset \mathbb R^n$ is said to be arc-symmetric if, for every real analytic arc $\gamma:(-1,1)\to \mathbb R^n$ such that $\gamma((-1,0))\subset X$, we have $\gamma((-1,1))\subset X$. Euclidean closed arc-symmetric subsets of $\mathbb{R}^n$ have nice topological and algebraic properties as they are the closed sets of a Noetherian topology of $\mathbb R^n$ finer than the usual Zariski topology of $\mathbb R^n$. Let $X\subset\mathbb R^n$ be a Euclidean closed arc-symmetric set, then a semialgebraic function $f:X\to \mathbb{R}$ is said to be arc-analytic if $f\circ\gamma:(-1,1)\to \mathbb{R}$ is real analytic for every real analytic arc $\gamma:(-1,1)\to X$. In \cite{BM90}, Bierstone and Milman showed that a semialgebraic function $f:X\to \mathbb{R}$ on a nonsingular real algebraic set $X$ of $\mathbb R^n$ is arc-analytic if and only if there is a finite sequence of blowings-up with nonsingular algebraic centers $\pi:\widetilde{X}\to X$ such that $f\circ\pi$ is a Nash function. Moreover, the vanishing locus of an arc-analytic function $f:\mathbb R^n\to \mathbb R$ is an arc-symmetric set, conversely it was recently shown that arc-symmetric sets are characterized as zero loci of arc-analytic functions, see \cite{AS17}. In addition, the ring of arc-analytic functions has good properties too as it satisfies the Nullstellensatz and any arc-analytic function on a Euclidean closed arc-symmetric set is the restriction of an arc-analytic function defined on the whole $\mathbb{R}^n$, see \cite{AS18} for the latter. However, although the topology is Noetherian, the ring of arc-analytic functions on a closed arc-symmetric set $X\subset \mathbb R^n$ is not Noetherian. \\
		
		More recently, the notion of \emph{$k$-regulous} function for $k\in\mathbb{N}\cup\{\infty\}$ has been introduced in an independent way, with slightly different definitions, by several authors with different motivations. Kucharz \cite{K09} introduced the notion of $\mathcal{C}^k$-rational map as follows: a $\mathcal C^k$ map $f:X\to Y$ between nonsingular irreducible real algebraic sets is said to be $\mathcal C^k$-rational if $f|_U$ is regular for some nonempty Zariski open set $U$ of $X$; $\mathcal C^\infty$-rational maps are shown to be regular on $X$. See also \cite{Ban24b,Kuc26} for recent progress. In the meantime, Kollár and Nowak \cite{KN15} introduced the notion of continuous hereditarily rational function, that is, a continuous function $f:X\to \mathbb{R}$ is  \emph{hereditarily rational} if every irreducible algebraic subset $Z$ of $X$ has a nonempty Zariski open subset $U\subset Z$ such that $f|_{U}$ is regular. The motivation for the latter definition is particularly subtle, one needs this hereditary property in order to define a well-behaved class of functions on singular algebraic sets. Moreover, the main result of the paper states that $f:X\to\mathbb R$ is a continuous hereditarily rational function if and only if there is a $\mathcal{C}^0$-rational (or regulous) function $F:\mathbb R^n\to \mathbb R$ such that $f=F|_X$. The notion of $k$-regulous function has then been introduced by Fichou, Huisman, Mangolte and the second author \cite{FHMM16} as follows, let $k\in\mathbb{N}\cup\{\infty\}$, a $\mathcal{C}^k$ function $f:\mathbb{R}^n\to \mathbb R$ is \emph{$k$-regulous} if there is a Zariski open dense subset $U$ of $\mathbb{R}^n$ such that $f|_{U}$ is regular. Observe that the latter definition coincides with $\mathcal{C}^k$-rational functions of $\mathbb R^n$ in the sense of \cite{K09}. In \cite{FHMM16} the authors study in depth the properties of the ring of $k$-regulous functions, both from an algebraic and a topological point of view, for $k\in \mathbb{N}$. They show that the zero sets of $k$-regulous functions on $\mathbb{R}^n$ coincide with the family of Euclidean closed Zariski-constructible subsets of $\mathbb{R}^n$ and that $k$-regulous functions are arc-analytic. Although the ring $\mathcal{R}^k(\mathbb{R}^n)$ of $k$-regulous functions on $\mathbb R^n$ is not Noetherian, the authors proved it has some very desirable properties: the Nullstellensatz holds for  $\mathcal{R}^k(\mathbb{R}^n)$ and it is radically principal. These two crucial properties allowed the authors to study the spectrum $\operatorname{Spec}\left(\mathcal{R}^k(\mathbb{R}^n)\right)$ together with its structure sheaf, which induces then a sheaf of $k$-regulous functions on $\mathbb{R}^n$ satisfying a $k$-regulous version of Cartan's Theorems A and B. A similar result is then obtained by restriction over any Euclidean closed Zariski-constructible subset of $\mathbb R^n$. It is worth mentioning that recently Banecki showed that a $k$-rational function $f:X\to \mathbb R$ in the sense of \cite{K09}, that is, defined on a nonsingular algebraic set $X\subset\mathbb R^n$, admits a $k$-regulous extension to the whole $\mathbb R^n$, see \cite{Ban25} for the case of surfaces and the preprint \cite{Ban24} for the general case. More recently, the definition of $k$-regulous functions has been also extended over any real closed field $R$ both in \cite{KK23} and \cite{BFM24}.  \\
		
		The aim of the present paper is to introduce a new class of real functions that we call \emph{Zariski-constructible functions of class $\mathcal{C}^k$} ($\ZC^k$-functions for short), where $k\in\mathbb N\cup\{\infty\}$, extending non-trivially the notion of $k$-regulous and regular function for any $k$, even in the case of functions on the real line and for $k=\infty$. Our functions are not arc-analytic in general but they are defined in a completely intrinsic way by asking that their graph be a Zariski-constructible set, in a similar manner as semialgebraic functions are real valued functions whose graph is semialgebraic. We will discuss in the next section more precisely our definitions and main results. As our techniques work over any real closed field $R$, we will directly define $\ZC^k$-functions over any real closed field $R$.
		
		\subsection{$\ZC$-sets, $\ZC^k$-functions and main results}
		
		Let us introduce the main notions and results of our paper. Let $R$ be a real closed field. We say that a set $X\subset R^n$ is a \emph{Zariski-constructible set} ($\ZC$-set for short) if it can be written as a Boolean combination of algebraic subsets of $R^n$. Equivalently, $X\subset R^n$ is a $\ZC$-set if and only if there exist $\ell\in \mathbb N$ and algebraic subsets $S_i$ and $T_i$ of $ R^n$ for $i\in\{1,\dots, \ell\}$ such that
		\[
		X=\bigcup_{i=1}^\ell (S_i^c\cap T_i).
		\]
		Observe that we can always assume the above union to be disjoint.
		
		There are two essential topological properties of $\ZC$-sets. The first is that Euclidean locally closed $\ZC$-sets are Euler, that is, the Euler characteristic of the link at any point is even, the second is the following result, which is crucial for our purposes.
		
		\begin{thm*}
			Euclidean closed $\ZC$-subsets of $ R^n$ are the closed sets of a Noetherian topology on $ R^n$ that we call the $\ZC$-topology of $R^n$.
		\end{thm*}
		
		In light of the previous result we will refer to Euclidean closed (resp. open) $\ZC$-subsets of $ R^n$ as $\ZC$-closed (resp. $\ZC$-open) subsets of $R^n$. The next example shows that the $\ZC$-topology is strictly finer than the Zariski topology of $R^n$.
		
		\begin{eg}\label{eg:Cartan}
			A Euclidean closed $\ZC$-set may not be algebraic. Indeed, let $$Z\coloneq \left\{(x,y,z)\in\mathbb R^3\,:\,z\left(x^2+y^2\right)=x^3\right\}$$ be Cartan's umbrella and let $X\coloneq Z\setminus\{(0,0,z)\in R^3\,:\,z\neq0\}$ be its canopy. Observe that $X\subset R^3$ is a $\ZC$-closed set but it is not an algebraic set, as Cartan's umbrella $Z\subset R^3$ is an irreducible algebraic set and $X$ contains the regular locus of $Z$, so the Zariski closure of $X$ coincides with $Z$.
		\end{eg}
		
		\begin{rem}\label{rem:C1}
			Let $C$ be an algebraically closed field of characteristic $0$. Observe that we can define in a similar way the notion of $\ZC$-subset $X$ of $C^n$. However, as the Euclidean closure of $X$ coincides with the Zariski closure of $X$, Euclidean closed $\ZC$-subsets of $C^n$ do coincide with algebraic subsets of $C^n$. This shows that the $\ZC$-topology of $R^n$ encodes intrinsically real phenomena, as in Example \ref{eg:Cartan}. 
		\end{rem}
		
		We point out that the above theorem  was already known in the case $R=\mathbb{R}$, see \cite{Par04}, but the proof uses Borel--Moore homology which is not directly available over a general real closed field $R$. Our argument is a completely combinatorial one relying on the triangulability of closed and bounded semialgebraic sets, which is available over any real closed field $R$. This is the main content of Subsection \ref{subsec:Euler} below.\\
			
		In Section \ref{sec:3}, we start by defining $\ZC$-maps and functions as follows. Let $X\subset R^n$ and $Y\subset R^m$. We say that $f:X\to Y$ is a \emph{Zariski-constructible map} ($\ZC$-map for short) if the graph $\Gamma_f\subset R^n\times R^m$ is a $\ZC$-set. Observe that, by definition, if $f:X\to Y$ is a bijective $\ZC$-map, then $f^{-1}:Y\to X$ is a $\ZC$-map too. When $Y=R$, we say that $f:X\to R$ is a $\ZC$-function. In Subsection \ref{sub:3.1} we prove additional properties of $\ZC$-maps in order to show that the set $\ZC(X)$ of $\ZC$-functions on $X$, endowed with the pointwise addition and multiplication, is a ring, see Proposition \ref{prop:ZC-ring}.\\
		
		In Subsection \ref{subsec:3.2}, we define the notion of $\ZC^k$-functions, first over $\ZC$-open subsets of $R^n$. Let $k\in\mathbb N\cup\{\infty\}$ and let $U\subset R^n$ be a $\ZC$-open set. We say that a function $f:U\to R$ is a \emph{Zariski-constructible function of class $\mathcal{C}^k$} (\emph{$\ZC^k$-function} for short) if it is a $\ZC$-function on $U$ of class $\mathcal C^k$. Denote by $\ZC^k(U)$ the subring of $\ZC(U)$ of $\ZC^k$-functions. Consider the following simple example showing the difference between $\ZC^0$-functions, regulous and continuous arc-analytic functions.
		
		\begin{eg}
			Consider the continuous function $f:R\to R$ defined as $f(x)=\sqrt[2r+1]{x}$. Observe that $f$ is a $\ZC$-function for every $r\in\mathbb N$, as its inverse is the polynomial function $p:R\to R$ defined as $p(x)=x^{2r+1}$ for every $r\in\mathbb N$. However, $f$ is neither regulous (as regulous functions on $R$ are regular) nor arc-analytic (when $R=\mathbb{R}$) for every $r\ge1$.
		\end{eg}
		
		More generally, we show that $k$-regulous functions on $R^n$ are $\ZC^k$-functions on $R^n$ for any $k\in\mathbb{N}\cup\{\infty\}$ (see Lemma \ref{lem:regulous} below), but that the latter class is strictly larger for any $k\in\mathbb{N}\cup\{\infty\}$. This is the case even for $k=\infty$, for which $\infty$-regulous functions on $R^n$ are just regular ones, and even for $n=1$, when $k$-regulous functions are just regular ones for any $k\in\mathbb{N}\cup \{\infty\}$. Indeed, Examples \ref{eg:i-throot}-\ref{eg:Nash} below show that we have the following strict inclusions:
		\begin{enumerate}
			\item $\ZC^\infty(R)\subsetneq\dots \subsetneq \ZC^{k+1}(R)\subsetneq\ZC^{k}(R)\subsetneq\dots\subsetneq\ZCc(R)$;
			\item $\mathcal{R}(R)\subsetneq \ZC^\infty(R)\subsetneq \mathcal{N}(R)$,
		\end{enumerate}
		where $\mathcal{N}(R)$ denotes the ring of Nash functions on $R$.\\
		
		Then, in Subsection \ref{subsec:ZCk(Rn)} we focus on the algebraic properties of the ring of $\ZC^k$-functions. We stress that such a ring have very different properties when $k\in\mathbb N$ compared to the case $k=\infty$. For the rest of the introduction we will assume $k\in\mathbb N$, however in Sections \ref{sec:3} and \ref{sec:4} some results are extended also for $k=\infty$. Specific algebraic properties of the ring of $\ZC^\infty$-functions are the subject of ongoing investigation by the authors.\\
		
		Let then $k\in\mathbb N$. We prove that $\ZC$-closed subsets of $R^n$ are exactly the zero sets of $\ZC^k$-functions on $R^n$ (see Proposition \ref{prop:zerolocus}). Although the ring $\ZC^k(R^n)$ is not Noetherian, see Example \ref{eg:notnoeth}, we obtain several good algebraic properties using the Noetherianity of the $\ZC$-topology together with a $\ZC^k$-version of the Łojasiewicz inequality (see Proposition \ref{prop:LojP}). Observe that Example \ref{eg:ki-root} below shows that Proposition \ref{prop:LojP} does not hold for $k=\infty$. The first main algebraic property we prove is the Nullstellensatz for the ring $\ZC^k(R^n)$. Let us introduce the necessary notation, similar to the usual one for zero sets and vanishing ideals in algebraic geometry.\\
		
		Let $k\in\mathbb{N}$. Let $S\subset \ZC^k(R^n)$ and $X\subset R^n$, we denote by
		\begin{align*}
			\mathcal{Z}(S)&\coloneq\{x\in R^n\,:\,f(x)=0\text{ for every $f\in S$}\},\\
			\I_{\ZC^k(R^n)}(X)&\coloneq\{f\in\ZC^k(R^n)\,:\,f(x)=0\text{ for every $x\in X$}\}
		\end{align*}
		the \emph{vanishing set} $\mathcal{Z}(S)$ of $S$ in $R^n$ and the \emph{vanishing ideal} $\I_{\ZC^k(R^n)}(X)$ of $X$ in $\ZC^k(R^n)$. When $S=\{f_1,\dots,f_s\}$ is finite, we shorten the notation as $\mathcal{Z}(f_1,\dots,f_s)\coloneq \mathcal{Z}(\{f_1,\dots,f_s\})$.\\
		
		\begin{thm*}[Nullstellensatz]
			Let $k\in\mathbb N$ and let $I$ be an ideal of $\ZC^k(R^n)$. Then
			\[
			\I_{\ZC^k(R^n)}\left(\mathcal Z(I)\right)=\sqrt{I}.
			\]
		\end{thm*}
		
		\noindent The Nullstellensatz has the following remarkable consequences:
		\begin{itemize}
			\item Let $I\subset \ZC^k(R^n)$ be an ideal. Then, $I$ is radical if and only if $I$ is real.
			\item The ring $\ZC^k(R^n)$ is \emph{radically principal}, that is, let $I\subset \ZC^k(R^n)$ be an ideal, then $\sqrt{I}=\sqrt{(f)}$ for some $f\in I$. Moreover, we show that $\sqrt{I}=\sqrt{(f)}$ if and only if $\mathcal{Z}(I)=\mathcal{Z}(f)$.
			\item Maximal ideals of $\ZC^k(R^n)$ are exactly the ideals of the form $\mathfrak{m}_a\coloneq\{f\in\ZC^k(R^n) : f(a)=0\}$ for some $a\in R^n$.
		\end{itemize}
		
		Then, we focus on the ring $\ZC^k(U)$, where $U\subset R^n$ is a $\ZC$-open set. As already mentioned, there is $f\in\ZC^k(R^n)$ such that $\mathcal Z(f)=R^n\setminus U$. We show that, for any such $f\in\ZC^k(R^n)$, restriction morphism $\rho:\ZC^k(R^n)\to \ZC^k(U)$ induces a canonical isomorphism $\rho_{f}:\ZC^k(R^n)_f\to \ZC^k(U)$ (see Proposition \ref{prop:localization}). Observe that Example \ref{eg:ki-root} shows that Proposition \ref{prop:localization} does not hold for $k=\infty$. As a consequence of the Nullstellensatz for $\ZC^k(R^n)$, Proposition \ref{prop:localization} and the usual properties of the localization, we get the Nullstellensatz also for the ring $\ZC^k(U)$ for any $\ZC$-open subset $U$ of $R^n$, see Corollary \ref{cor:Nullstellensatz-U}. \\
		
		Then we define $\ZC^k$ functions on $\ZC$-locally closed sets. A $\ZC$-locally closed subset of $R^n$ is a $\ZC$-subset of $R^n$ which is locally closed with respect to the $\ZC$-topology of $R^n$. Let $X\subset R^n$ be a $\ZC$-locally closed set. We say that $f:X\to R$ is a \emph{Zariski-constructible function of class $\mathcal{C}^k$} (\emph{$\ZC^k$-function} for short) if it is locally, for the $\ZC$-topology, the restriction of a $\ZC^k$-function defined on a $\ZC$-open subset of $R^n$, see Definition \ref{def:ZC^k(X)}. We denote by $\ZC^k(X)$ the ring of $\ZC^k$-functions on $X$. Observe that if $f:X\to R$ is a $\ZC^k$-function, then $f$ is a $\ZC$-function and it is continuous with respect to the induced Euclidean topology on $X$ (see Remark \ref{rem:ZC-continuity}).\\
		
		\begin{rem}\label{rem:C2}
			Let us have a brief comparison with the complex setting. Let $X\subset \mathbb C^n$ be an algebraic set and $f:X\to \mathbb C$ be a continuous $\ZC$-function, that is, a continuous function (with respect to the Euclidean topology) whose graph $\Gamma_f\subset \mathbb{C}^n\times \mathbb C$ is a $\ZC$-set. As $f:X\to \mathbb C$ is continuous, then $\Gamma_f\subset \mathbb{C}^n\times \mathbb C$ is a $\ZC$-closed set, and hence an algebraic set by Remark \ref{rem:C1}. Therefore, a continuous $\ZC$-function $f:X\to \mathbb C$ is actually a continuous rational function as shown in \cite{Be26}.
		\end{rem}
		
		Let $X\subset R^n$ be a $\ZC$-locally closed set and denote by $\zcl{X}$ the closure of $X$ in the $\ZC$-topology of $R^n$. As a consequence of the Nullstellensatz for the ring $\ZC^k(U)$, where $U\subset R^n$ is a $\ZC$-open subset of $R^n$ such that $X=U\cap\zcl{X}$, we have the following extension theorem for $\ZC^k$-functions on $X$.
		
		\begin{thm*}[Extension of $\ZC^k$-functions]
			Let $k\in\mathbb N$. Let $X\subset R^n$ be a $\ZC$-locally closed set and let $U\subset R^n$ be a $\ZC$-open set such that $X=U\cap\zcl{X}$. Let $f\in\ZC^k(X)$, then there exists $F\in\ZC^k(U)$ such that $F|_{X}=f$.
		\end{thm*}
		
		Observe that, if $X\subset R^n$ is a $\ZC$-closed set we may choose $U=R^n$ in the statement of the above theorem.\\
		
		Then, in Section \ref{sec:4} we study the sheaf of $\ZC^k$-functions on a $\ZC$-locally closed set $X\subset R^n$. More precisely, let $X\subset R^n$ be a $\ZC$-locally closed set and endow $X$ with the induced $\ZC$-topology. We define the sheaf of rings $\ZC^k_X$ on $X$ as follows: for every $\ZC$-open subset $U$  of $X$ we define $\ZC^k_X(U)\coloneq \ZC^k(U)$, see Lemma \ref{lem:sheaf}.
		
		We show that the $\ZC$-closed subsets of $X$ are in correspondence with the Zariski-closed subsets of $\operatorname{Spec}\left(\ZC^k(X)\right)$ (see Theorem \ref{thm:correspondence} and Corollary \ref{cor:correspondence}). Denote then by $\jmath_{X}:X\hookrightarrow \operatorname{Spec}\left(\ZC^k(X)\right)$ the inclusion map such that $\jmath_{X}(a)$ is the maximal ideal of $\ZC^k$-functions on $X$ vanishing at $a\in X$. We show that $\jmath_{X}(X)=\operatorname{Spec}_m\left(\ZC^k(X)\right)$, the set of maximal ideals of $\ZC^k(X)$. Then, we show that the sheaf $\ZC^k_X$ on $X$ is isomorphic to the sheaf $\jmath_{X}^{-1}\left(\mathcal{O}_{\widetilde{X}}\right)$, where $\mathcal{O}_{\widetilde{X}}$ denotes the usual structure sheaf of $\operatorname{Spec}\left(\ZC^k(X)\right)$, see Corollary \ref{cor:inverse-sheaf}.\\
		
		Therefore, an application of \cite[Thm.\,1.3.1]{Gro61} gives the following $\ZC^k$-version of Cartan's Theorems A and B in the algebraic formulation of Serre (see \cite[\S45,Thm.\,2;\,\S46,Thm.\,3]{Ser55}). We refer the reader to Subsection \ref{sub:Cartan} for the precise definitions involved in the statement.
		
		\begin{thm*}[Cartan's Theorems A and B]
			Let $k\in\mathbb{N}$. Let $X\subset R^n$ be a $\ZC$-locally closed set and let $\mathcal{F}$ be a quasi-coherent $\ZC^k_X$-module on $X$. Then:
			\begin{enumerate}[label=(\Alph*),nosep]
				\item $\mathcal{F}$ is generated by its global sections.
				\item $H^p(X,\mathcal{F})=0$ for every $p>0$.
			\end{enumerate}
		\end{thm*}
		
		Finally, we conclude the paper by proposing an alternative proof of the extension theorem relying just on Cartan's Theorem B for quasi-coherent $\ZC^k_X$-modules on $R^n$.

		\section{Topological properties of $\ZC$-sets}
		
		\subsection{Main results}
		
		The main goal of this section is to study the topological properties of $\ZC$-subsets of $R^n$. Our main result in this section, which is of crucial importance for the rest of the paper, is the following:
		
		\begin{thm}\label{thm:ZC-top}
			Euclidean closed $\ZC$-subsets of $ R^n$ are the closed sets of a Noetherian topology on $ R^n$ that we call the $\ZC$-topology of $R^n$.
		\end{thm}
		
		The proof of Theorem \ref{thm:ZC-top} is postponed until later in the section, as it relies on technical properties of Euler sets introduced and studied in Subsection \ref{subsec:Euler} below. Let us first deduce some basic properties of the $\ZC$-topology of $R^n$ after introducing the necessary notation.\\

		Let $X\subset R^n$. We denote by $\ecl{X}$ the \emph{Euclidean closure} of $X$ in $R^n$. Observe that, as $X$ is in particular semialgebraic, then $\ecl{X}$ is a semialgebraic subset of $R^n$ by \cite[Prop.\,2.2.2]{BCR}. Denote by $\zcl{X}\subset R^n$ and by $\zarcl{X}\subset R^n$ the \emph{$\ZC$-closure} and the \emph{Zariski closure} of $X$ in $R^n$, respectively. As the Zariski topology of $R^n$ is coarser than the $\ZC$-topology of $R^n$, which is coarser than the Euclidean topology of $R^n$, we have that $X\subset\ecl{X}\subset\zcl{X}\subset\zarcl{X}$. Recall that the \emph{dimension} $\dim(X)$ of a semialgebraic set $X\subset R^n$ is the dimension of the Zariski closure $\zarcl{X}$ as an algebraic set; equivalently, it is  the largest $d\in\mathbb N$ such that $X$ contains a Nash submanifold of $R^n$ of dimension $d$, see \cite[\S2.8]{BCR}.
		
		\begin{prop}\label{prop:dim}
			Let $X\subset R^n$ be a semialgebraic set. Then $\dim\left(\zcl{X}\right)=\dim (X)$.
		\end{prop}
		\begin{proof}
			As $X\subset\zcl{X}\subset\zarcl{X}$ we have that $\zarcl{X}=\zarcl{\zcl{X}}$. Therefore, $\dim\zcl{X}=\dim X$.
		\end{proof}
		
Let $X\subset R^n$. Denote by $\ebd{X}\coloneq\ecl{X}\setminus\eint{X}$ and $\zbd{X}\coloneq \zcl{X}\setminus\zint{X}$ the \emph{Euclidean boundary} of $X$ and the \emph{$\ZC$-boundary} of $X$, respectively. Observe that, if $X\subset R^n$ is a $\ZC$-set, then $\ebd X$ is semialgebraic and $\zbd{X}$ is a $\ZC$-set.
\begin{prop}\label{prop:closure}
Let $X\subset R^n$ be a $\ZC$-set. Then
\[
\zcl{X}=X\cup\zcl{\ecl{X}\setminus X}=X\cup\zcl{\ebd{X}}.
\]
Besides, if $X$ is nonempty then $\dim\left(\zcl{X}\setminus X\right)<\dim X$.
\end{prop}
\begin{proof}
Set $F\coloneq\ecl{X}\setminus X$ (which is semialgebraic). Since $\ecl{X}=X\cup F$, we have that $X\cup\zcl{F}=\ecl{X}\cup\zcl{F}$ is a $\ZC$-closed set containing $X$, so that $\zcl{X}\subset X\cup\zcl{F}$. Moreover, $F\subset\ebd{X}\subset\ecl{X}\subset\zcl{X}$ so that
\[
\zcl{X}\subset X\cup\zcl{F}\subset X\cup\zcl{\ebd{X}}\subset\zcl{X}.
\]
This shows that $\zcl{X}=X\cup\zcl{F}=X\cup\zcl{\ebd{X}}$.

In particular, $\zcl{X}\setminus X\subset\zcl{F}$. Therefore, by Proposition \ref{prop:dim},
\[
\dim\left(\zcl{X}\setminus X\right)\le\dim\left(\zcl{F}\right)=\dim(F)<\dim(X),
\]
where the last inequality comes from \cite[Prop.\,2.8.13]{BCR}.
\end{proof}

\begin{rem}\label{rem:whitney}
For a $\ZC$-set $X$, the inclusion $\ecl{X}\setminus X\subset\zcl{X}\setminus X$ can be strict. Consider for instance the Whitney umbrella without the $z$-axis $X\coloneq\{(x,y,z)\in R^3\,:\,zx^2=y^2,\,x^2+y^2\neq 0\}$. Observe that $\ecl{X}\setminus X=\{(x,y,z)\in R^3\,:\,x^2+y^2=0,\,z\geq 0\}$ and that
\[
\zcl{X}\setminus X=\zcl{\ecl{X}\setminus X}=L\supsetneq\ecl{X}\setminus X,
\]
where $L$ is the $z$-axis. In particular, this shows that $\ecl{X}$ is a Euclidean closed semialgebraic set that is not a $\ZC$-set as it is strictly contained in  $\zcl{X}=X\cup L=\zarcl{X}$ that is, the Whitney umbrella.
\end{rem}

		\subsection{A detour to Euler sets and their combinatorial properties}\label{subsec:Euler}
		
		Let us settle the basic combinatorial definitions we will use later on in the section. Let $a_0,\dots,a_d$ be $d+1$ affinely independent points of $R^n$. A \emph{$d$-simplex} $\sigma\coloneq[a_0,\dots,a_d]$ of $R^n$ with vertices $a_0,\dots,a_d$ is the set of $a\in R^n$ such that there are non-negative $\lambda_0,\dots,\lambda_d\in R$ such that $\sum_{i=0}^{d} \lambda_i=1$ and $a=\sum_{i=0}^{d} \lambda_i a_i$. We say that a $d$-simplex has \emph{dimension} $d$. Let $\{a_{i_0},\dots,a_{i_k}\}\subset\{a_0,\dots,a_d\}$, then the simplex $\tau\coloneq[a_{i_0},\dots,a_{i_k}]\subset\sigma$ is a $k$-\emph{face} of $\sigma$. A \emph{finite simplicial complex} of $R^n$ is a finite collection $K=(\sigma_i)_{i=1,\dots,\ell}$ such that the faces of each simplex $\sigma_i\in K$ belong to $K$ and if $\sigma_i,\sigma_j\in K$, then $\sigma_{i}\cap\sigma_j$ is either empty or a common face of $\sigma_i$ and $\sigma_j$. { A subset $K'$ of a finite simplicial complex $K$ is called a \emph{simplicial subcomplex} of $K$ if it is a (finite) simplicial complex of $R^n$ as well.} Denote by $|K|\coloneq\bigcup_{i=1}^\ell \sigma_i\subset R^n$ the \emph{realization} of the simplicial complex $K$. Denote by $K^{(d)}$ the set of simplices of $K$ of dimension $d$. {To simplify the notation, we will identify $\sigma\in K^{(0)}$ with the unique $a_0\in R^n$ such that $\sigma=[a_0]$.} \\
		
		{ Let $\sigma,\tau\subset R^n$ be two simplices of $R^n$ such that $\sigma\cap\tau=\varnothing$, where the empty set is intended as a simplex of $R^n$ of dimension $-1$. Define the \emph{join} $\sigma\ast\tau$ of $\sigma$ and $\tau$ as the $(\dim(\sigma)+\dim(\tau)+1)$-simplex spanned by the vertices of $\sigma$ and $\tau$. In particular, $\sigma\ast\varnothing=\sigma$ and $\varnothing\ast \tau=\tau$. Let $K$ be a finite simplicial complex. Let $S$ and $T$ be (finite) simplicial subcomplexes of $K$ such that $S\cap T=\varnothing$. Define the \emph{join} $S\ast T$ of $S$ and $T$ in $K$ as the (finite) simplicial subcomplex of $K$:
			\[
			S\ast T\coloneq\{\sigma\ast \tau\in K : \sigma\in S, \tau\in T\},
			\]
		where the empty set is admitted as a simplex both of $S$ and $T$.\\
		
		Let $K$ be a finite simplicial complex of $R^n$ and $\sigma\in K$. We define the \emph{star} $\st_\sigma(K)$ of $\sigma$ in $K$ as $\st_\sigma(K)\coloneq\{\tau\in K : \sigma \subset \tau\}$. Define the \emph{boundary} $\dst_\sigma(K)$ of the star $\st_\sigma(K)$ of $\sigma$ in $K$ as
		\[
		\dst_\sigma(K)\coloneq\{\tau\in K\,:\,\sigma\nsubset\tau\text{ and }\exists\rho\in K,\,\tau,\sigma\subset\rho\}.
		\]
		Define the \emph{link} of $\sigma$ in $K$ as 
		\[
		\lk_\sigma(K)\coloneq\{\tau\in K\,:\,\sigma\cap\tau=\varnothing\text{ and }\exists\rho\in K,\,\tau,\sigma\subset\rho\}.
		\]
		Observe that both $\dst_\sigma(K)$ and $\lk_\sigma(K)$ are (finite) simplicial subcomplexes of $K$. If $\dim\sigma=0$ then $\dst_\sigma(K)=\lk_\sigma(K)$, otherwise, if $\dim\sigma>0$, then $\partial\sigma\cap\lk_\sigma(K)=\varnothing$ and
		\[
		\dst_\sigma(K)=\partial\sigma\ast\lk_\sigma(K), 
		\]
		where $\partial \sigma$ denotes the (finite) simplicial subcomplex of $K$ constituted by the $k$-faces of $\sigma$ for every $k<\dim\sigma$.} \\
		
		Recall that every closed and bounded semialgebraic set $S\subset R^n$ admits a finite semialgebraic triangulation \cite[Thm.9.2.1]{BCR}, that is, a finite simplicial complex $K\coloneq(\sigma_i)_{i=1,\dots,\ell}$ and a semialgebraic homeomorphism $\psi:|K|\to S$. Therefore, we may define the homology groups $H_k(S,\mathbb{Z}_2)$ of a closed and bounded semialgebraic set $S$ of $R^n$ as the simplicial homology $H_k(K,\mathbb{Z}_2)$. For more information on simplicial homology and for the independence of $H_k(S,\mathbb{Z}_2)$ from the choice of the triangulation $\psi:|K|\to S$ we refer to \cite{Mun84} and \cite[\S 11.7]{BCR}.\\
		
		Let $S\subset R^n$ be a locally closed semialgebraic set. Recall that the \emph{link} of $S$ at $x\in S$ is defined as $\lk_x(S)\coloneq S\cap \mathbb S^n_\varepsilon(x)$ for $0<\varepsilon<<1$. As the local structure of a semialgebraic set at any point is conic \cite[Thm.9.3.6]{BCR}, we obtain that the link $\lk_x(S)$ of $S$ at $x$ is well-defined up to semialgebraic homeomorphism. Moreover, as $S$ is locally closed, $\lk_x(S)\coloneq S\cap \mathbb S^{n-1}_\varepsilon(x)$ is a closed semialgebraic subset of $R^n$ for $\varepsilon$ sufficiently small.\\
		
		We are now in position to introduce the main object of this subsection.
		
		\begin{defn}\label{def:euler}
			Let $S\subset R^n$ be a locally closed semialgebraic set. We say that $S$ is \emph{Euler} if
			\[
			\chi(\lk_x(S))\coloneq\sum_i(-1)^i\dim_{\mathbb Z_2}\left(H_i(\lk_x(S),\mathbb Z_2)\right)
			\]
			is even for every $x\in S$.
		\end{defn}
		
		{ In 1971 Sullivan \cite{Sul71} proved that every real algebraic subset $X$ of $\mathbb{R}^n$ is Euler. Let us sketch the argument. Let $X$ be an algebraic subset of $\mathbb R^n$ and denote by $X_{\mathbb C}$ its complexification in $\mathbb C^n$. Then, $\chi(\lk_x(X))\equiv\chi(\lk_x(X_{\mathbb C}))\mod 2$ since the fixed points of the conjugation on $\lk_x(X_{\mathbb C})$ coincide with $\lk_x(X)$. Finally, $\chi(\lk_x(X_{\mathbb C}))=0$ as $\lk_x(X_{\mathbb C})$ admits a stratification made of odd dimensional strata only. Later on, this argument has been extended over any real closed field, see \cite[Thm.\,11.2.2]{BCR}.}
		
		\begin{lemma}\label{lem:ZCEuler}
			Euclidean locally closed $\ZC$-sets are Euler.
		\end{lemma}
		\begin{proof}
			Let $X\subset R^n$ be a locally closed $\ZC$-set. Let $x\in X$ then, for $\varepsilon>0$ small enough, the link $\lk_x(X)$ is well defined and semialgebraically homeomorphic to $X\cap\mathbb{S}_\varepsilon(x)$, which is a closed subset of $R^n$. Write $X=\bigsqcup_{i=1}^r S_i^c\cap T_i$ where $S_i,T_i$ are real algebraic sets, then $X\cap\mathbb{S}_\varepsilon(x)=\bigsqcup_{i=1}^r (S_i^c\cap T_i)\cap\mathbb S_{\varepsilon}(x)$. Therefore 
			\begin{align*}
				\chi(\lk_x(X))&=\chi(X\cap \mathbb S_\varepsilon(x))=\sum_{i=1}^r\chi((S_i^c\cap T_i)\cap\mathbb S_{\varepsilon}(x))\\&=\sum_{i=1}^r\chi(\lk_x(T_i))-\chi(\lk_x(S_i))\equiv0\mod2,
			\end{align*}
			as $T_i$ and $S_i$ are algebraic subsets of $R^n$ for every $i\in\{1,\dots,r\}$, see \cite[Thm.\,11.2.2]{BCR}. This shows that $X$ is Euler, as desired.
		\end{proof}
		
		Let $S\subset R^n$ be a locally closed semialgebraic set and let $x\in S$. Let $\dot{S}\subset R^m$ be the semialgebraic Alexandrov compactification of $S$, with $m\geq n$, see \cite[Prop.2.5.9 \& Def.2.5.11]{BCR}. Denote by $\iota:S\hookrightarrow\dot{S}$ the semialgebraic inclusion map of $S$ in $\dot{S}$ and denote by $\{\infty\}\coloneq\dot{S}\setminus\iota(S)$. As above, by \cite[Thm.9.2.1]{BCR}, there exists a triangulation of $\dot{S}$ compatible with $\{\infty\}$, that is, a finite simplicial complex $K\coloneq(\sigma_i)_{i=0,1,\dots,\ell}$ and a semialgebraic homeomorphism $\psi:|K|\to \dot{S}$ such that $\psi(\sigma_{0})=\infty$.\\
		
		In what follows, we characterize the property of being Euler for a locally closed semialgebraic set $S\subset R^n$ with combinatorial properties of its Alexandrov compactification $\dot{S}\subset R^m$, for $m\geq n$.
		
		\begin{lemma}\label{lem:parity}
			Let $S\subset R^n$ be a closed and bounded semialgebraic set. Let $K\coloneq(\sigma_i)_{i=1,\dots,\ell}$ be a finite simplicial complex of $R^n$ and $\psi:|K|\to S$ be a semialgebraic triangulation of $S$. Then, $S$ is Euler if and only if $\#\{\tau\in K\,:\,\sigma\subsetneq\tau\}$ is even for each simplex $\sigma\in K$. 
		\end{lemma}
		
		\begin{proof}
			Let $\sigma\in K$ and $x\in\relint(\sigma)$. Consider the finite simplicial subcomplexes $\st_{\sigma}(K)\coloneq\{\tau\in K : \sigma\subset\tau\}$ and $\dst_\sigma(K)\coloneq\{\tau\in K\,:\,\sigma\nsubset\tau\text{ and }\exists\rho\in K,\,\tau,\sigma\subset\rho\}$. Consider the finite simplicial complex $K'$ of $R^n$, sometimes called the star subdivision of $K$ at $x$, defined as
			\[
			K'\coloneq(K\setminus\st_{\sigma}(K))\cup\partial\st_{\sigma}(K)\cup\{x\}\cup\{x\ast \tau\subset R^n\,:\,\tau\in\partial\st_{\sigma}(K)\}.
			\] 
			Observe that if $\dim(\sigma)=0$, then $K=K'$. Moreover, by definition of $K'$, we have that $|K|=|K'|$ and $\dst_x(K')=\dst_\sigma(K)$ (see Figure \ref{fig:link}). As the homeomorphism $\psi:|K'|\to S$ induces the local conic structure of $S$ at $\psi(x)$, we directly get that $\chi(\lk_{\psi(x)}(S))=\chi(\dst_\sigma(K))$. Therefore, as every $x\in|K|$ is contained in $\relint(\sigma)$ for a unique $\sigma\in K$, we showed that $S$ is Euler if and only if $\chi(\dst_\sigma(K))$ is even for every simplex $\sigma\in K$.\\
			
			\begin{figure}
				\centering
				\tikzset{every picture/.style={line width=0.75pt}}
				
				\begin{tikzpicture}[x=0.5pt,y=0.5pt,yscale=-1,xscale=1]
					
					\draw [color={rgb, 255:red, 0; green, 0; blue, 0 }  ,draw opacity=1 ][line width=1.5]    (228,122) -- (485,214) ;
					\draw [shift={(356.5,168)}, rotate = 19.7] [color={rgb, 255:red, 0; green, 0; blue, 0 }  ,draw opacity=1 ][fill={rgb, 255:red, 0; green, 0; blue, 0 }  ,fill opacity=1 ][line width=1.5]      (0, 0) circle [x radius= 4.36, y radius= 4.36]   ;
					\draw [line width=1.5]    (228,122) -- (185,25) ;
					\draw [shift={(185,25)}, rotate = 246.09] [color={rgb, 255:red, 0; green, 0; blue, 0 }  ][fill={rgb, 255:red, 0; green, 0; blue, 0 }  ][line width=1.5]      (0, 0) circle [x radius= 4.36, y radius= 4.36]   ;
					\draw [shift={(228,122)}, rotate = 246.09] [color={rgb, 255:red, 0; green, 0; blue, 0 }  ][fill={rgb, 255:red, 0; green, 0; blue, 0 }  ][line width=1.5]      (0, 0) circle [x radius= 4.36, y radius= 4.36]   ;
					\draw [line width=1.5]    (185,25) -- (485,214) ;
					\draw [line width=1.5]  [dash pattern={on 1.69pt off 2.76pt}]  (228,122) -- (435,101) ;
					\draw [line width=1.5]    (435,101) -- (485,214) ;
					\draw [shift={(485,214)}, rotate = 66.13] [color={rgb, 255:red, 0; green, 0; blue, 0 }  ][fill={rgb, 255:red, 0; green, 0; blue, 0 }  ][line width=1.5]      (0, 0) circle [x radius= 4.36, y radius= 4.36]   ;
					\draw [shift={(435,101)}, rotate = 66.13] [color={rgb, 255:red, 0; green, 0; blue, 0 }  ][fill={rgb, 255:red, 0; green, 0; blue, 0 }  ][line width=1.5]      (0, 0) circle [x radius= 4.36, y radius= 4.36]   ;
					\draw [line width=1.5]    (228,122) -- (175.29,190.2) ;
					\draw [line width=1.5]    (175.29,190.2) -- (485,214) ;
					\draw [shift={(485,214)}, rotate = 4.39] [color={rgb, 255:red, 0; green, 0; blue, 0 }  ][fill={rgb, 255:red, 0; green, 0; blue, 0 }  ][line width=1.5]      (0, 0) circle [x radius= 4.36, y radius= 4.36]   ;
					\draw [shift={(175.29,190.2)}, rotate = 4.39] [color={rgb, 255:red, 0; green, 0; blue, 0 }  ][fill={rgb, 255:red, 0; green, 0; blue, 0 }  ][line width=1.5]      (0, 0) circle [x radius= 4.36, y radius= 4.36]   ;
					\draw [line width=1.5]  [dash pattern={on 1.69pt off 2.76pt}]  (228,122) -- (368.73,262.18) ;
					\draw [line width=1.5]    (368.73,262.18) -- (485,214) ;
					\draw [shift={(368.73,262.18)}, rotate = 337.49] [color={rgb, 255:red, 0; green, 0; blue, 0 }  ][fill={rgb, 255:red, 0; green, 0; blue, 0 }  ][line width=1.5]      (0, 0) circle [x radius= 4.36, y radius= 4.36]   ;
					\draw [line width=1.5]    (435,101) -- (326,112) ;
					\draw [line width=1.5]    (306.73,200.18) -- (368.73,262.18) ;
					\draw    (175.29,190.2) -- (356.5,168) ;
					\draw [color={rgb, 255:red, 0; green, 0; blue, 0 }  ,draw opacity=1 ] [dash pattern={on 0.84pt off 2.51pt}]  (356.5,168) -- (368.73,262.18) ;
					\draw [color={rgb, 255:red, 0; green, 0; blue, 0 }  ,draw opacity=1 ]   (185,25) -- (356.5,168) ;
					\draw  [dash pattern={on 0.84pt off 2.51pt}]  (435,101) -- (356.5,168) ;
					\draw    (361.73,207.18) -- (368.73,262.18) ;
					\draw    (435,101) -- (381,147) ;
					
					\draw (375,183.4) node [anchor=north west][inner sep=0.75pt]    {$\sigma $};
					\draw (348,142.4) node [anchor=north west][inner sep=0.75pt]    {$x$};
					\draw (384,56.4) node [anchor=north west][inner sep=0.75pt]    {$|\st_{\sigma}(K)|$};
					\draw (105,135) node [anchor=north west][inner sep=0.75pt]    {$|\st_{x}(K')|$};

				\end{tikzpicture}
				\caption{Comparison between $\st_{\sigma}(K)$ and $\st_{x}(K')$.}
				\label{fig:link}
			\end{figure}
			
			{Let $\sigma\in K$. As $\dst_\sigma(K)=\partial\sigma\ast\lk_\sigma(K)$, we have that $\dst_\sigma(K)^{(k)}=\partial\sigma^{(k)}\sqcup\lk_\sigma(K)^{(k)}\sqcup L^{(k)}$, where $L^{(k)}\coloneq \{\sigma'\ast\tau'\,:\,\sigma'\in\partial\sigma^{(s)}, \tau'\in \lk_\sigma(K)^{(t)},\,s+t+1=k,\,s,t\geq 0\}$ for every $k\in \mathbb{N}$, and hence:
			\begin{align*}
				\chi(\dst_\sigma(K))&\coloneq\sum_{k}(-1)^k \# \left(\partial\sigma^{(k)}\sqcup\lk_\sigma(K)^{(k)}\sqcup L^{(k)}\right)\\
													&=\chi(\partial\sigma)+\chi(\lk_\sigma(K))+ \sum_{k} (-1)^k \bigg(\sum_{s+t+1=k} \#\partial\sigma^{(s)}\cdot \#\lk_\sigma(K)^{(t)}\bigg) \\
													&=\chi(\partial\sigma)+\chi(\lk_\sigma(K))-\chi(\partial\sigma)\chi(\lk_\sigma(K)).
			\end{align*}}
		
			As $\sigma$ is a simplex, we have that $\chi(\partial\sigma)=0$ if $\dim(\sigma)$ is even and $\chi(\partial\sigma)=2$ if $\dim(\sigma)$ is odd, and therefore:
			\[
			\chi(\dst_\sigma(K))\equiv \chi(\lk_\sigma(K)) \mod 2.
			\]
			Observe that, by definition of $\lk_\sigma(K)$, there is a bijection between the $k$-simplices of $\lk_\sigma(K)$ and the $(\dim(\sigma)+k+1)$-simplices of $K$ containing $\sigma$, and therefore:
			\begin{align*}
				\chi(\lk_\sigma(K))&=\sum_{k}(-1)^{k}\#\{\tau\in K\,:\,\tau\in K^{(\dim(\sigma)+k+1)},\,\sigma\subset\tau\} \\
				&\equiv\#\{\tau\in K\,:\,\sigma\subsetneq\tau\}\mod2.
			\end{align*}
			This shows that $S$ is Euler if and only if $\#\{\tau\in K\,:\,\sigma\subsetneq\tau\}$ is even for every $\sigma\in K$, as desired.
		\end{proof}

		Let us relate the property of being Euler for a locally closed semialgebraic set $S\subset R^n$ with its semialgebraic Alexandrov compactification $\dot{S}\subset R^m$ with $m\geq n$.
		
		\begin{lemma}\label{lem:PtEuler}
			Let $S$ be a Euclidean locally closed semialgebraic set. Then $S$ is Euler if and only if its semialgebraic Alexandrov compactification $\dot{S}$ is Euler.
		\end{lemma}
		\begin{proof}
			Assume that $\dot{S}$ is Euler and let $\iota:S\hookrightarrow\dot{S}$ be the inclusion map, which is a semialgebraic homeomorphism when restricted to its image $\dot{S}\setminus\{\infty\}$. Let $x\in S$ and $0<\varepsilon<<1$ such that $\lk_{\iota(x)}(\dot{S})\coloneq\dot{S}\cap \mathbb S_\varepsilon(\iota(x))\subset \iota(S)$, then $\chi(\lk_x(S))=\chi(\lk_{\iota(x)}(\iota(S)))=\chi(\lk_{\iota(x)}(\dot{S}))\equiv0\mod2$. Thus $S$ is Euler.
			
			Conversely, assume that $S$ is Euler. Let $\psi:|K|\to\dot{S}$ be a finite triangulation compatible with $\{\infty\}\coloneq\dot{S}\setminus\iota(S)$ and denote by $\sigma_0$ the $0$-simplex such that $\psi(\sigma_0)=\infty$. For every $\sigma\in K$, set $n_\sigma\coloneq\#\{\rho\,:\,\rho\subsetneq\sigma\}$ and $m_\sigma\coloneq\#\{\tau\,:\,\sigma\subsetneq\tau\}$. Then
			\[
			\sum_{\sigma\in K}m_\sigma=\sum_{\sigma\in K}n_\sigma=\sum_{\sigma\in K}(2^{\dim\sigma+1}-2)\equiv0\mod2.
			\]
			Let $\sigma\neq\sigma_{0}$ and $x\in \relint(\sigma)$. As shown in the proof of Lemma \ref{lem:parity}, $\chi(\lk_{\psi(x)}(\dot{S}))$ is even if and only if $m_\sigma$ is even. As $\lk_{\psi(x)}(\dot{S})=\lk_{\psi(x)}({\iota(S)})$ and $S$ is Euler, we obtain that $m_\sigma$ is even for every $\sigma\in K$ with $\sigma\neq\sigma_0$. Therefore, $m_{\sigma_0}=\sum_{\sigma\in K}m_\sigma-\sum_{\sigma\neq\sigma_0}m_\sigma$ is even as well. Finally, the application of Lemma \ref{lem:parity} ensures that $\dot{S}$ is Euler, as desired.
		\end{proof}

		\subsection{Proof of Theorem \ref{thm:ZC-top}}
		The general strategy of the proof presented below is adapted from \cite[Thm.\,1.4]{Kur88} which studies the case of semialgebraic arc-symmetric sets. Nonetheless, the first step (see \cite[\S1.6]{Kur88}), needs to be totally reworked for the $\ZC$-setting, it is the purpose of Lemma \ref{lem:top-manifold} below.
		
		For alternative proofs of Theorem \ref{thm:ZC-top} in the special case $R=\mathbb R$, using different methods, we refer to \cite[Prop.\,1.5]{Kur99}, \cite[Thm.\,2.8]{Par04}, and \cite[Prop.\,3.2]{KZ18}. Here we present an alternative combinatorial argument working over any real closed field.\\
		
		{ Let us recall a useful notion in real algebraic geometry over general real closed fields. Let $S\subset R^n$ be a semialgebraic set and endow $S$ with the induced Euclidean topology. We say that $S$ is \emph{semialgebraically connected} if $S_1$ and $S_2$ are disjoint closed subsets of $S$ such that $S=S_1\cup S_2$, then either $S_1=\varnothing$ or $S_2=\varnothing$. Moreover, a semialgebraic set $S\subset R^n$ is semialgebraically connected if and only if for every $x,y\in S$ there is a semialgebraic continuous path $\gamma:[0,1]\to X$ such that $\gamma(0)=x$ and $\gamma(1)=y$, see \cite[Prop.\,2.5.13]{BCR}.}
			
		\begin{rem}
			{When $R=\mathbb R$, a semialgebraic set $S\subset R^n$ is semialgebraically connected if and only if $S$ is connected in the usual sense, see \cite[Thm.\,2.4.5]{BCR}. However, when considering the real closed field $R=\overline{\mathbb Q}^r\coloneq \overline{\mathbb{Q}}\cap \mathbb{R}$, that is, the field of real algebraic numbers, we have that $R$ is semialgebraically connected but totally disconnected, see \cite[Ex.\,2.4.1(a)]{BCR}.}
		\end{rem}
		
		{Let $S\subset R^n$ be a semialgebraic set of dimension $d$. We say that $S$ is a \emph{semialgebraic topological manifold} if for every $x\in S$ there is a semialgebraic open neighborhood $U$ of $x$ in $S$ and a semialgebraic homeomorphism $\varphi:\mathbb{B}^d\to U$, where $\mathbb{B}^d\subset R^d$ denotes the standard unit ball.}
		
		\begin{lemma}\label{lem:top-manifold}
			Let $X\subset R^n$ be a semialgebraically connected semialgebraic topological manifold and $Y\subset X$ be a semialgebraic set closed in $X$. If $Y$ is Euler then either $Y= X$ or $\dim Y < \dim X$.
		\end{lemma}
		
		The above Lemma \ref{lem:top-manifold} is proven in \cite[Lem.\,1.6]{Par04}, in the case $R=\mathbb{R}$, using Borel--Moore homology. Below, we present an alternative proof relying solely on the triangulation of semialgebraic sets, which is available over any real closed field.
		
		\begin{proof}

			Assume that $\dim X=\dim Y=d$ and that $X$ is unbounded (if $X$ is bounded the following proof suitably simplified works without semialgebraic Alexandrov compactification).
			
			Let $\dot{X}\coloneq \iota(X)\cup\{\infty\}$ be the semialgebraic one-point compactification of $X$. Then $\dot{X}$ and $\dot{Y}\coloneq \iota(Y)\cup\{\infty\}$ are semialgebraic, $\dot{Y}$ is a closed subset of $\dot{X}$, since $Y$ is closed in $X$, and $\dot{Y}$ is Euler by Lemma \ref{lem:PtEuler}.
			
			Let $K\coloneq(\sigma_i)_{i=0,\dots,\ell}$ be a finite simplicial complex and let $\psi:|K|\to \dot{X}$ be a finite triangulation of $\dot{X}$ compatible with $\dot{Y}$ and $\{\infty\}$. Denote by $L$ the simplicial subcomplex of $K$ such that $\psi(|L|)=\dot{Y}$ and $\sigma_0\in K$ such that $\psi(\sigma_0)=\infty$. For $e\in\mathbb N$, denote by $K^{(e)}$ and by $L^{(e)}$ the set of simplices of dimension $e$ of $K$ and $L$, respectively. Note that $L^{(d)}\neq\varnothing$ since $\dim X=\dim Y=d$.
			
			\begin{claim}
				For every $\sigma\in L^{(d)}$ and $\tau\in K^{(d)}$, if $\sigma\cap\tau\in L^{(d-1)}$ then $\tau\in L^{(d)}$.
			\end{claim}
			Let  $\sigma\in L^{(d)}$ and $\tau\in K^{(d)}$ be such that $\sigma'\coloneq\sigma\cap\tau\in L^{(d-1)}$. Observe that $\#\{\rho\in L^{(d)}\,:\,\sigma'\subset\rho\}$ is even and strictly positive, as $Y$ is Euler and $\sigma\in\{\rho\in L^{(d)}\,:\,\sigma'\subset\rho\}$. Besides, since $X$ is a semialgebraic topological manifold and $\dot{X}$ is its one-point compactification, there are exactly two simplices $\sigma_1,\sigma_2\in K^{(d)}$ such that $\sigma'\subset\sigma_1\cap\sigma_2$. Therefore, $\{\sigma,\tau\}=\{\rho\in L^{(d)}\,:\,\sigma'\subset\rho\}$, which proves the claim. \cqed
			
			\begin{claim}
				$K^{(d)}=L^{(d)}$.
			\end{claim}
			{ Let $\sigma\in L^{(d)}$ and let $\tau\in K^{(d)}\setminus L^{(d)}$. Pick $x\in\relint(\sigma)$ and $y\in\relint(\tau)$. Observe that $\psi(x),\psi(y)\in\dot{X}\setminus\{\infty\}$. As $\dot{X}\setminus\{\infty\}=\iota(X)$ is semialgebraically connected, then $|K|\setminus\{\sigma_0\}=\psi^{-1}(\dot{X}\setminus\{\infty\})$ is semialgebraically connected as well, therefore there is a continuous semialgebraic path $\gamma:[0,1]\to |K|\setminus\{\sigma_0\}$ such that $\gamma(0)=x$ and $\gamma(1)=y$. Consider then a sequence of simplices $\tau_1,\dots,\tau_s\in K^{(d)}$ such that $\tau_1\coloneq \sigma$, $\tau_s\coloneq \tau$, $\tau_i\cap\tau_{i+1}\in K^{(d-1)}$ and $\gamma([0,1])\subset \tau_1\cup \ldots\cup\tau_s$.
			The above claim shows that $\tau_i\in L^{(d)}$ for every $i\in\{1,\dots,s\}$. Therefore $\tau=\tau_s\in L^{(d)}$ and the claim is proven.\cqed \\}
			
			Finally, since $X$ is a semialgebraic topological manifold, every simplex of $K^{(e)}$ for $e\leq d-1$ is contained in some simplex of dimension $d$, therefore $\dot{X}\setminus\{\infty\}=\psi(|K|)\setminus\{\infty\}=\psi(|L|)\setminus\{\infty\}=\dot{Y}\setminus\{\infty\}$. This shows that $Y=X$.
		\end{proof}
		
		\begin{lemma}\label{lem:tmind}
			Let $V,W\subset R^n$ be algebraic sets and let $X$ be a semialgebraically connected component of $V\setminus W$. Let $(Y_i)_{i\in I}$ be a family of Euclidean closed $\ZC$-subsets of $R^n$. Then there exists $J\subset I$ finite such that
			\[
			X\cap\bigcap_{i\in I}Y_i=X\cap\bigcap_{i\in J}Y_i.
			\]
		\end{lemma}
		\begin{proof}
			Let us prove the result by induction on $d\in\mathbb N$. If $d=\dim X=0$ there is nothing to prove. Let $d>0$ and assume the statement holds when $\dim X<d$.
				
			Let $X\subset R^n$ such that $d=\dim X$. First of all, let us make a reduction. Let $\zarcl{X}$ denote the Zariski closure of $X$ in $R^n$. Then, we may suppose $V=\zarcl{X}$ and therefore $d\coloneq\dim X=\dim \zarcl{X}$. Moreover, as $V$ is algebraic, up to substituting $Y_i$ with $Y_i\cap V$, we may suppose $Y_i\subset V$ for every $i\in I$.
			
			If $X\subset Y_i$ for all $i\in I$, then there is nothing to prove. So we may assume that $X\nsubset Y_{i_0}$ for some $i_0\in I$. By the previous reduction we have that $Y_i\subset V=\zarcl{X}$ for every $i\in I$, hence Lemma \ref{lem:top-manifold} gives that $\dim X\cap Y_{i_0}<d$. As $X\cap Y_{i_0}$ is semialgebraic, \cite[Prop.9.1.8]{BCR} ensures the existence of a finite Nash stratification $\{X_j\}_{j=1}^s$ of $X\cap Y_{i_0}$, that is, each $X_j$ is a semialgebraically connected Nash submanifold of $ R^n$ such that $X\cap Y_{i_0}=\bigcup_{j=1}^{s} X_j$ and $X_{j_1}\cap X_{j_2}=\varnothing$ for every $j_1,j_2\in\{1,\dots,s\}$ with $j_1\neq j_2$.

			Then
			\[
			X\cap\bigcap_{i\in I} Y_i=X\cap Y_{i_0}\cap\bigcap_{i\in I\setminus\{i_0\}} Y_i=\bigcup_{j=1}^s\left(X_j\cap\bigcap_{i\in I\setminus\{i_0\}} Y_i\right).
			\]
			As $X_j$ is a semialgebraically connected Nash manifold such that $\dim X_j<d$ for every $j\in\{1,\dots,s\}$, by the inductive assumption, for every $j\in\{1,\dots,s\}$ there is a finite set $I_j\subset  I\setminus\{i_0\}$ such that 
			\[
			X_j\cap\bigcap_{i\in I\setminus\{i_0\}} Y_i=X_j\cap\bigcap_{i\in I_j}Y_i.
			\]
			Set $J_0\coloneq\bigcup_{j=1}^s I_j$ and $J\coloneq J_0\cup\{i_0\}$, then $J_0$ and $J$ are finite and
			\[
			X\cap\bigcap_{i\in I} Y_i=\bigcup_{j=1}^s\Bigg(X_j\cap\bigcap_{i\in I_j}Y_i\Bigg)=\bigcap_{i\in J_0}\left(Y_i\cap\bigcup_{j=1}^sX_j\right)=\bigcap_{i\in J_0}\left(Y_i\cap X\cap Y_{i_0}\right)=X\cap\bigcap_{i\in J}Y_i,
			\]
			as desired.
		\end{proof}
		
	\begin{proof}[Proof of Theorem \ref{thm:ZC-top}]
		Let $\{Y_i\}_{i\in I}$ be a family of Euclidean closed $\ZC$-subsets of $ R^n$. By Lemma \ref{lem:ZCEuler} and Lemma \ref{lem:tmind} (with $X= R^n$), there exists $J\subset I$ finite such that $\bigcap_{i\in I}Y_i=\bigcap_{i\in J}Y_i$. This shows that Euclidean closed $\ZC$-subsets of $R^n$ are closed under intersection, hence they are the closed sets of a topology on $ R^n$. { Moreover, as we showed that any arbitrary intersection of Euclidean closed $\ZC$-subsets of $R^n$ reduces to a finite one, we actually proved that the $\ZC$-topology is Noetherian.}
	\end{proof}
		
		\begin{rem}\label{rem:irrcomps}
			 As the $\ZC$-topology is Noetherian, we have the notions of irreducible $\ZC$-closed set, the irreducible components and the irreducible decomposition of a $\ZC$-closed set $X$ of $R^n$.
		\end{rem}
		
		\begin{notation}
			In light of Theorem \ref{thm:ZC-top} we fix the following notation: we will refer to a Euclidean closed $\ZC$-subset and open $\ZC$-subset of $R^n$ as a \emph{$\ZC$-closed} and \emph{$\ZC$-open} subset of $R^n$, respectively.
		\end{notation}

		\section{$\ZC^k$-maps and functions}\label{sec:3}
		
		\subsection{$\ZC$-maps and functions: fundamental properties}\label{sub:3.1}
		
		The notion of $\ZC$-set provides a natural notion of $\ZC$-function and $\ZC$-map, in a similar way as the notion of semialgebraic set provides a natural notion of semialgebraic function and map.
		
		\begin{defn}[$\ZC$-map \& $\ZC$-function]
			Let $X\subset R^n$ and $f:X\to R^m$ be a map. We say that $f$ is a \emph{$\ZC$-map} if the graph $\Gamma_f\subset R^n\times R^m$ is a $\ZC$-set. If $m=1$, we say that $f$ is a \emph{$\ZC$-function}. We denote by $\ZC(X,R^m)$ the set of $\ZC$-maps from $X$ to $R^m$ and by $\ZC(X)\coloneq\ZC(X,R)$ the set of $\ZC$-functions on $X$.
		\end{defn}
			
		The next lemma may be thought as a weak analogue of Tarski--Seidenberg (or of Chevalley's theorem in the complex case) for $\ZC$-mappings.
		
		{
		\begin{lemma}\label{lem:injective}
			Let $X\subset R^n\times R^m$ be a $\ZC$-set. Denote by $\pi:R^n\times R^m\to R^n$ the projection onto the first $n$ coordinates. If $\pi|_{X}$ is injective, then $\pi(X)\subset R^n$ is a $\ZC$-set.
		\end{lemma}
		\begin{proof}
		Let $X=\bigcup_{i=1}^r S_i^c\cap T_i$ where $S_i,T_i$ are real algebraic subsets of $R^n\times R^m$. As $\pi|_X$ is injective, we have that $\pi(X)=\bigcup_{i=1}^r \pi(S_i^c\cap T_i)$, therefore it suffices to prove the statement for $X=S^c\cap T$, where $S$ and $T$ are real algebraic subsets of $R^n\times R^m$. Let $s\in R[x_1,\dots,x_n,y_1,\dots,y_m]$ be a polynomial such that $s^{-1}(0)=S$. Consider the real algebraic subset
		\[
		Y\coloneq\{(x,y,z)\in T\times R : zs(x,y)=1\}
		\]
		of $R^n\times R^m\times R$ and the regular map $\widetilde{\pi}:R^n\times R^m\times R\to R^n$ defined as $\widetilde{\pi}\coloneq \pi\circ \pi_{n+m}$, where $\pi_{n+m}:R^n\times R^m\times R\to R^n\times R^m$ is the projection $\pi_{n+m}(x,y,z)=(x,y)$. Observe that, by definition of $Y$, the map $\pi_{n+m}|_Y$ is injective, and therefore $\widetilde{\pi}|_{Y}$ is so. This reduces to prove the statement for an algebraic set $X\subset R^n\times R^m$.
		
		\vspace{0.5em}
		
		Assume $X\subset R^n\times R^m$ is an algebraic set. Dealing separately with each irreducible component of $X$, we may suppose $X$ is irreducible. Let us show that $\pi(X)\subset R^n$ is a $\ZC$-set by induction on $d\coloneq\dim(X)$. The case $d=0$ is evident, therefore suppose $d>0$. Denote by $V\coloneq\zarcl{\pi(X)}\subset R^n$ the Zariski closure of $\pi(X)$ in $R^n$. Observe that $V$ is irreducible, as $\pi$ is a polynomial map and $X$ is so. As $\pi|_X$ is injective and the algebraic dimension of $V$ coincides with the semialgebraic dimension of $\pi(X)$, we have that $\dim(V)=d$.
		
		\vspace{0.5em}
		
		Denote by $C\coloneq R[i]$ the algebraic closure of $R$. Consider the complexifications $X_C\subset C^{n+m}$ and $V_C\subset C^n$ of $X$ and $V$, respectively. Recall that $X_C$ and $V_C$ are irreducible, as $X$ and $V$ are so, and $d=\dim(V_C)=\dim(X_C)$. Denote by $\pi_C:C^n\times C^m\to C^n$ the projection onto the first $n$ coordinates. Let us show first that $\pi_C|_{X_C}:X_C\to V_C$ is dominant, that is, that $\overline{\pi_C(X_C)}^{C-\text{Zar}}=V_C$, where $\overline{\pi_C(X_C)}^{C-\text{Zar}}$ denotes the Zariski closure of $\pi_C(X_C)$ in $C^n$. Let $Z\coloneq \overline{\pi_C(X_C)}^{C-\text{Zar}}$ and denote by $\sigma:C^n\to C^n$ the complex conjugation $\sigma(x_1,\dots,x_n)=(\overline{x}_1,\dots,\overline{x}_n)$. As $V_C$ is the complexification of the real algebraic set $V\subset R^n$, we have that $Z'\coloneq Z\cup \sigma(Z)\subset V_C$ is an algebraic subset of $C^n$ defined over $R$. Therefore, $\zarcl{\pi(X)}\subset Z'\subset V_C$ and hence $V_C=(\zarcl{\pi(X)})_C\subset Z'\subset V_C$. This shows that $V_C=Z'\coloneq Z\cup\sigma(Z)$ and hence $Z\coloneq \overline{\pi_C(X_C)}^{C-\text{Zar}}=V_C$ as $V_C$ is irreducible.
		
		\vspace{0.5em}
	
		Consider the induced field homomorphism $\theta_\pi:C(V_C) \to C(X_C)$ defined as $\theta_\pi(f)\coloneq f\circ \pi$, where $C(V_C)$ and $C(X_C)$ denote the field of rational functions on $V_C$ and $X_C$, respectively. As $\pi_C:X_C\to V_C$ is dominant and $\dim(X_C)=d=\dim(V_C)$, we have that $\theta_\pi$ is injective, the field extension $C(X_C)|C(V_C)$ is algebraic and its degree $e\coloneq[C(X_C)|C(V_C)]$ coincides with the cardinality of the generic fiber of $\pi_C$ by \cite[Prop.7.16]{Har92}.
		
		\vspace{0.5em}
			
		Denote by $W\subset V_C$ an algebraic subset of $C^n$ such that $\#(\pi_C^{-1}(y)\cap X_C)=e$ for every $y\in V_C\setminus W$. Observe that $\dim(W)<d$ as $V_C$ is irreducible and hence $V'\coloneq W\cap R^n\subset V$ is an algebraic subset of $R^n$ of dimension $<d$. This shows that $U\coloneq V\setminus V'$ is a nonempty Zariski open subset of $V$. Observe that $\pi_C^{-1}(y)\cap (X_C\setminus X)$ consists of pairs of complex conjugated points for every $y\in V$, then $\# (\pi^{-1}(y)\cap X)\equiv e \mod 2$ for every $y\in U$. Therefore, as $\pi|_X$ is injective and $V\coloneq\zarcl{\pi(X)}$, we obtain that $e$ is odd and $U\subset \pi(X)$, as $\# (\pi^{-1}(y)\cap X)\in\{0,1\}$ for every $y\in V$. Consider the algebraic subset $X'\coloneq \pi^{-1}(V')\cap X$ of $R^n$. As $\pi|_{X}$ is injective and $X'\subset X$, then $\pi|_{X'}$ is injective. Moreover, as $\dim(V')<d$, we have that $\dim(X')<d=\dim(X)$ and hence $\pi(X')$ is a $\ZC$-subset of $R^n$ by inductive assumption. Therefore:
		\[
		\pi(X)=\pi(X)\cap V=\pi(X)\cap (U\sqcup V')=(\pi(X)\cap U)\sqcup (\pi(X)\cap V')=U\sqcup \pi(X')
		\]
		is a $\ZC$-subset of $R^n$, as desired.
		\end{proof}

		The proposition below, in the case $R=\mathbb R$, is proven in \cite[Thm.\,2.8]{Par04} using algebraically constructible functions in the sense of \cite{PS97,MP97}. We provide an alternative proof working over any real closed field $R$ relying just on Lemma \ref{lem:injective}. 
			
		\begin{prop}\label{prop:ZCf}
			Let $f:X\to R^m$ be a $\ZC$-map, where $X\subset R^n$. Then:
			\begin{enumerate}[label=(\roman*),nosep]
				\item\label{ZCf-3} The domain $X$ of $f$ is a $\ZC$-set.
				\item\label{ZCf-2} If $Z\subset R^m$ is a $\ZC$-set then $f^{-1}(Z)$ is too.
				\item\label{ZCf-1} If $f$ is injective and if $Y\subset X$ is a $\ZC$-set then $f(Y)$ is too.
			\end{enumerate}
		\end{prop}
		
		\begin{proof}
		Denote by $\Gamma_f\subset R^n\times R^m$ the graph of $f$.
		
		\vspace{0.5em}
		
		\ref{ZCf-3} As $f$ is a map, the projection $\pi:R^n\times R^m\to R^m$ onto the first $n$ coordinates is injective when restricted to $\Gamma_f$, and therefore $X=\pi(\Gamma_f)$ is a $\ZC$-subset of $R^n$ by Lemma \ref{lem:injective}.
		
		\vspace{0.5em}
		
		\ref{ZCf-2} Let $Z\subset R^m$ be a $\ZC$-set. Consider the map $f|_{f^{-1}(Z)}:f^{-1}(Z)\to R^m$. Observe that the graph $\Gamma_{f|_{f^{-1}(Z)}}$ of $f|_{f^{-1}(Z)}$ coincides with $\Gamma_f\cap (R^n\times Z)$, which is a $\ZC$-subset of $R^n\times R^m$. Therefore, $f|_{f^{-1}(Z)}$ is a $\ZC$-map and then $f^{-1}(Z)$ is a $\ZC$-subset of $R^n$ by \ref{ZCf-3}.
		
		\vspace{0.5em}
		
		\ref{ZCf-1} Assume $f:X\to R^m$ is an injective $\ZC$-map and let $Y\subset X$ a $\ZC$-subset of $R^n$. Denote by $g:f(X)\to X$ the inverse of $f$. Observe that $\Gamma_g\subset R^m\times R^n$ is a $\ZC$-set as $\Gamma_f\subset R^n\times R^m$ is so, and therefore $f(Y)=g^{-1}(Y)$ is a $\ZC$-subset of $R^m$ by \ref{ZCf-2}.
		\end{proof}
	}
		
		Let us describe some basic properties of $\ZC$-maps and functions.
		
		\begin{prop}\label{prop:composition}
			Let $f:X\to R^m$ and $g:Y\to R^p$ be two $\ZC$-maps such that $f(X)\subset Y$. Then $g\circ f:X\to R^p$ is a $\ZC$-map.
		\end{prop}
		\begin{proof}
			Let $f:X\to R^m$ and $g:Y\to R^p$ be two $\ZC$-maps with $X\subset R^n$, $Y\subset R^m$ and $f(X)\subset Y$. Denote by $\pi:R^n\times R^m\times R^p\to R^n\times R^p$ the projection map $\pi(x,y,z)\coloneq(x,z)$. Then
			\[
			\Gamma_{g\circ f}=\pi\left(\left(\Gamma_f\times R^p\right)\cap\left(R^n\times\Gamma_g\right)\right).
			\]
			Observe that $\left(\Gamma_f\times R^p\right)\cap\left(R^n\times\Gamma_g\right)$ is a $\ZC$-subset of $R^n\times R^m\times R^p$. Moreover, $\pi|_{\left(\Gamma_f\times R^p\right)\cap\left(R^n\times\Gamma_g\right)}$ is injective as $(x,y,z)\in\left(\Gamma_f\times R^p\right)\cap\left(R^n\times\Gamma_g\right)$ if and only if $(x,y,z)=(x,f(x),g(f(x)))$, and therefore $\Gamma_{g\circ f}$ is a $\ZC$-set by Lemma \ref{lem:injective}. This shows that $g\circ f\in\ZC(X)$, as desired.\\
		\end{proof}
		
		\begin{lemma}\label{lem:product}
			Let $X\subset R^n$ and $Y\subset R^m$ be non-empty $\ZC$-sets. Let $f:X\to R^s$ and $g:Y\to R^t$ and denote by $f\times g:X\times Y\to R^s\times R^t$ the map $(f\times g)(x,y)\coloneq(f(x),g(y))$. Then, $f\in\ZC(X,R^s)$ and $g\in\ZC(Y,R^t)$ if and only if $f\times g\in\ZC(X\times Y,R^s\times R^t)$.
		\end{lemma}
		\begin{proof}
			Assume first that $f\times g\in\ZC(X\times Y,R^s\times R^t)$ and fix $x_0\in X$ and $y_0\in Y$. Then, $f(x)=\pi_s((f\times g)(\iota(x)))$ and $g(y)=\pi_t((f\times g)(\jmath(y)))$ where $\pi_s:R^s\times R^t\to R^s$ and $\pi_t:R^s\times R^t\to R^t$ denote the projections onto the first and second factor, respectively, and $\iota:X\to X\times Y$ and $\jmath:Y\to X\times Y$ denote the polynomial maps $\iota(x)\coloneq(x,y_0)$ and $\jmath(y)\coloneq(x_0,y)$ for every $x\in X$ and $y\in Y$. Therefore, $f\in\ZC(X,R^s)$ and $g\in\ZC(Y,R^t)$ by Proposition \ref{prop:composition}.\\
			
			On the other way round, assume $f\in\ZC(X,R^s)$ and $g\in\ZC(Y,R^t)$, that is, $\Gamma_f\subset R^n\times R^s$ and $\Gamma_g\subset R^m\times R^t$ are $\ZC$-sets. Then, $\Gamma_{f\times g}=s(\Gamma_f \times \Gamma_g)\subset (R^n\times R^m)\times(R^s\times R^t)$, where $\Gamma_f \times \Gamma_g\subset (R^n\times R^s)\times(R^m\times R^t)$ is a $\ZC$-set and $s:(R^n\times R^s)\times(R^m\times R^t)\to(R^n\times R^m)\times(R^s\times R^t)$ is the injective polynomial map $s(x,u,y,w)\coloneq(x,y,u,w)$. Therefore, $\Gamma_{f\times g}\subset (R^n\times R^m)\times(R^s\times R^t)$ is a $\ZC$-set by Proposition \ref{prop:ZCf}\ref{ZCf-1}.
		\end{proof}
		
		\begin{prop}\label{prop:ZC-ring}
			Let $X\subset R^n$ be a non-empty $\ZC$-set. The set of $\ZC$-functions $\ZC(X)$ on $X$ endowed with the pointwise addition and multiplication is a ring. 
		\end{prop}
		\begin{proof}
			Let $f,g\in\ZC(X)$. Let $\iota:R\to R\times R$ be the polynomial map $\iota(x)\coloneq(x,x)$. Then, $f+g=a\circ(f\times g)\circ\iota$ and $fg=m\circ(f\times g)\circ\iota$ where $a: R\times R\to R$ and $m: R\times R\to R$ are the polynomial maps defined by $a(x,y)=x+y$ and $m(x,y)=xy$. Therefore $f+g,fg\in\ZC(X)$ by Proposition \ref{prop:composition} and Lemma \ref{lem:product}.
		\end{proof}

		\subsection{$\ZC^k$-maps and functions: definition and first properties}\label{subsec:3.2}
		
		In this section, we add differentiability conditions to $\ZC$-functions. In order to define $\ZC$-functions of class $\mathcal{C}^k$, we first restrict ourselves to functions defined on $\ZC$-open subsets of $R^n$; $\ZC^k$-functions on $\ZC$-locally closed subsets of $R^n$ will be introduced in Definition \ref{def:ZC^k(X)}.
		
		\begin{defn}[$\ZC^k$-map \& $\ZC^k$-function]\label{def:ZCk}
			Let $U\subset R^n$ be a $\ZC$-open set and $f:U\to R^m$ be a map. We say that $f\in\ZC(U,R^m)$ is a \emph{$\ZC^k$-map} for $k\in\mathbb{N}\cup\{\infty\}$ if $f$ is a $\ZC$-map of class $\mathcal{C}^k$. If $m=1$, we say that such a $f$ is a \emph{$\ZC^k$-function}. We denote by $\ZC^k(U,R^m)$ the set of $\ZC^k$-maps from $U$ to $R^m$ and by $\ZC^k(U)\coloneq\ZC^k(U,R)$ the set of $\ZC^k$-functions on $U$.
		\end{defn}
		
		\begin{rem}\label{rmk:zeros}
			{Let $U\subset R^n$ be a $\ZC$-open set. Endow $U$ with the induced topology from the $\ZC$-topology of $R^n$. Let $g\in\ZCc(U)$ and let $Y\subset R$ be a $\ZC$-closed subset of $R$. Then, $g^{-1}(Y)$ is a $\ZC$-closed subset of $U$ by Proposition \ref{prop:ZCf}\ref{ZCf-2} and by continuity of $g$ with respect to the Euclidean topologies of $U$ and $R$. Therefore, $g$ is continuous with respect to the $\ZC$-topologies of  $U$ and $R$. In particular, if $Y=\{0\}$, then $\mathcal{Z}(g)\coloneq g^{-1}(0)$ is a $\ZC$-closed subset of $U$.}
		\end{rem}
		
		As a direct consequence of Proposition \ref{prop:composition}\ and Lemma \ref{lem:product} we have a ring structure also in $\ZC^k(U)$ for every $\ZC$-open subset $U$ of $R^n$ and for every $k\in\mathbb{N}\cup\{\infty\}$. 
		
		\begin{cor}\label{cor:ZCk-ring}
			Let $U\subset R^n$ be a $\ZC$-open set. Then the set of $\ZC^k$-functions $\ZC^k(U)$ on $U$ is a subring of $\ZC(U)$ for every $k\in\mathbb{N}\cup\{\infty\}$.
		\end{cor}
		
		Denote by $\mathcal{R}(R)$ the usual ring of regular functions on $R$ and by $\mathcal{N}(R)$ the ring of Nash functions on $R$, that is, the ring of semialgebraic functions on $R$ of class $\mathcal{C}^\infty$. The next examples show that even in the simplest case of $U=R$ the ring $\ZC^k(R)$ strictly extends the ring $\mathcal{R}(R)$ of regular functions on $R$ for every $k\in\mathbb{N}\cup\{\infty\}$. 
		
		\begin{eg}\label{eg:i-throot}
			Consider the function $f:R\to R$ defined as $f(x)\coloneq\sqrt[2r+1]{x}$ for $r\in\mathbb{N}\setminus\{0\}$. Observe that its inverse $f^{-1}:R\to R$ is the polynomial map $f^{-1}(x)\coloneq x^{2r+1}$, therefore $\Gamma_f\subset R\times R$ is a $\ZC$-set and hence $f\in\ZCc(R)$. Evidently, $f$ is not of class $\mathcal{C}^1$ and hence $f\in\ZCc(R)\setminus\ZC^1(R)$.
		\end{eg}
		
		\begin{eg}\label{eg:ki-root}
			Fix $k,r\in\mathbb{N}\setminus\{0\}$ and consider the function $f_k:R\to R$ defined as $f_k(x)\coloneq x^k\sqrt[2r+1]{x}$. Observe that $f_k\in\ZCc(R)$ by Corollary \ref{cor:ZCk-ring} and Example \ref{eg:i-throot}. Moreover, as $f_k$ is of class $\mathcal{C}^k$ but not of class $\mathcal{C}^{k+1}$, we deduce that $f_k\in\ZC^k(R)\setminus\ZC^{k+1}(R)$. Observe that such an example arises as a special case of the $\ZC^k$-Łojasiewicz inequality, see Proposition \ref{prop:LojP} below.
		\end{eg}
		
		\begin{eg}\label{eg:ZCinf}
			Consider the algebraic set $Y\coloneq \{(x,y)\in R^2 : p(x,y)=y^3+2x^2y-x^4=0\}$. Consider the polynomial $p(x,y)\in R[x][y]$. Observe that the discriminant $\Delta_p(x)$ with respect to the variable $y$ is $\Delta_p(x)=-32x^6-27x^8<0$ for every $x\neq 0$ and $p(0,y)=0$ if and only if $y=0$. Therefore, for every $x\in R$ there is a unique root $(x,y)\in R^2$ of the polynomial $p(x,y)$ and, by Cardano's formulas, the $\mathcal{C}^0$ function $f:R\to R$ defined as
			\[
			f(x)=\sqrt[3]{\frac{x^4}{2}+\sqrt{\frac{x^8}{4}+\frac{8x^6}{27}}}+\sqrt[3]{\frac{x^4}{2}-\sqrt{\frac{x^8}{4}+\frac{8x^6}{27}}}
			\]
			satisfies $\Gamma_f=Y$. Let us show that actually $f\in\mathcal{C}^\infty(R)$ but $f\notin\mathcal{R}(R)$, and therefore $f\in\ZC^\infty(R)\setminus\mathcal{R}(R)$.
			
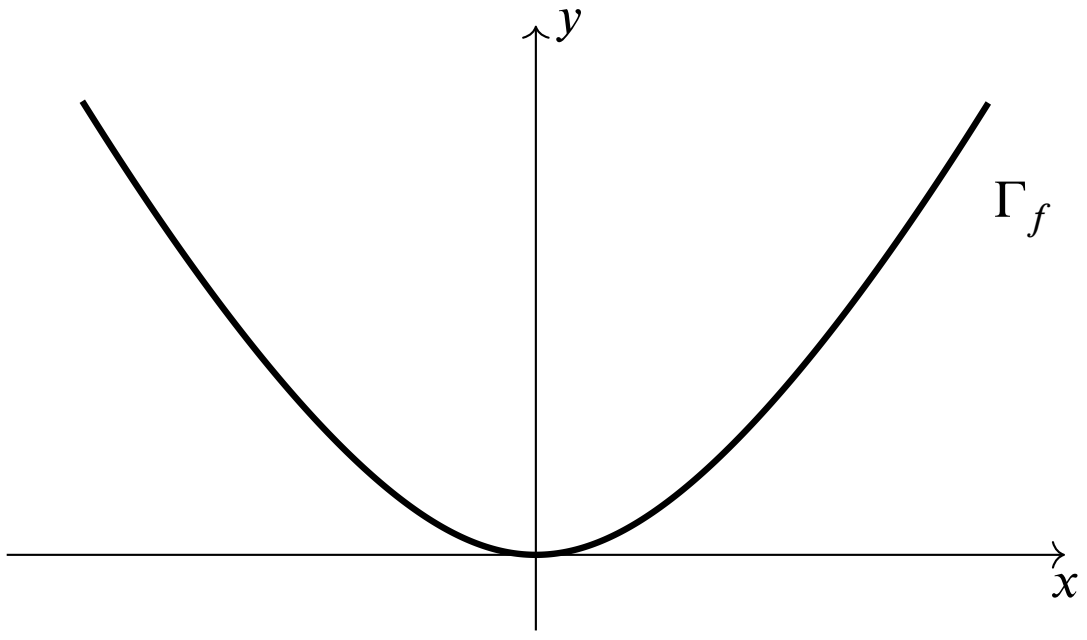
			\begin{center}
				\begin{tikzpicture}[x=1cm, y=1cm]\label{fig:ZC^infty}
					\draw[->] (-3.5,0) -- (3.5,0) node[below] {$x$};
					\draw[->] (0,-0.5) -- (0,3.5) node[right] {$y$};
					\draw[very thick] plot[domain=-1:1, samples=100, smooth] ({(\x)^3+2*\x}, {(\x)^4+2*(\x)^2});
					\draw (3,2.5) node [anchor=north west][inner sep=0.75pt]    {$\Gamma_f$};
				\end{tikzpicture}
			\end{center}
			
			Let us first show that $f\in\mathcal{C}^\infty(R)$. Of course, the only point to check is $x=0$. Observe that, by substituting $t=x^2$ and resolving the polynomial equation $-y^3-2yt+t^2=0$ with respect to the variable $t$ we get that $t=y\pm |y|\sqrt{1+y}$. Moreover, the unique solution in $t$ corresponding to the set $Y=\Gamma_f$ is $t=y(1+\sqrt{1+y})$ with $y\geq 0$. Consider the function $\phi:R_{\geq 0}\to R_{\geq 0}$ defined as $\phi(y)=y(1+\sqrt{1+y})$. Observe that $\phi$ admits a Nash extension to a neighborhood $U$ of $0$ in $R$ (as soon as $1+y>0$) and $\phi(0)=0$. In addition, $\phi'(y)=1+\sqrt{1+y}+\frac{y}{2\sqrt{1+y}}$ and $\phi'(0)=2\neq 0$, therefore, up to shrinking $U$ around $0$ if necessary, the Inverse Function Theorem \cite[Thm.2.9.7]{BCR} provides a neighborhood $V$ of $\phi(0)=0$ in $R$ and a Nash function $\psi:V\to U$ such that $\psi(\phi(y))=y$ for every $y\in U$. As $\phi(y)=t=x^2$, we get that $f(x)=y=\psi(\phi(y))=\psi(x^2)$ is a Nash function in a neighborhood of $0$.
			
			Let us show that $f\notin\mathcal{R}(R)$. Consider the polynomial $p(x,y)\coloneq y^3+2x^2y-x^4$ as an element in $R[x][y]$ and let us show that $p(x,y)$ does not have any regular (even rational) root $\alpha(x)\in\mathcal{R}(R)$. Let $A\coloneq R[x]$ and $K\coloneq \textnormal{Frac}(A)=R(x)$. Observe that $p$ is a monic polynomial in $A[y]$ and the ring $A$ is an UFD, therefore an application of the Integral Root Test \cite[\S IV,Prop.3.3]{Lan02} gives that if $\alpha(x)\in K$ satisfies $p(x,\alpha(x))=0$, then $\alpha(x)\in A\coloneq R[x]$. To show that $p(x,y)$ does not have any polynomial root $\alpha(x)\in R[x]$ it suffices to argue on the degree of $p(x,\alpha(x))=\alpha(x)^3+2x^2\alpha(x)-x^4$. Observe that $\alpha(x)=0$ is not a root, so denote by $d\coloneq \deg(\alpha)$. If $d=1$, then $\deg(p(x,\alpha(x)))=4$ and whenever $d\geq 2$ we have that $\deg(p(x,\alpha(x)))=3d$. Therefore, $p(x,\alpha(x))\neq 0$ for every $\alpha(x)\in R[x]$. This completes the proof that the above function $f\in\ZC^\infty(R)$ is not regular.
		\end{eg}
		
		\begin{eg}\label{eg:Nash}
			Consider the function $f:R\to R$ defined as $f(x)\coloneq \sqrt{1+x^2}$. Observe that $f$ is of class $\mathcal{C}^\infty$ and $\Gamma_f=\{(x,y)\in R\times R\,:\,y^2-x^2=1, y>0\}$ is a semialgebraic set, therefore $f$ is a Nash function on $R$. However, $f\notin \ZC^\infty(R)$ as $\Gamma_f$ is one of the two branches of the hyperbola $\{(x,y)\in R\times R\,:\,y^2-x^2=1\}$, which is an irreducible $\ZC$-set.
		\end{eg}
		
		As already mentioned, the previous Examples \ref{eg:i-throot},\ref{eg:ki-root},\ref{eg:ZCinf}\&\ref{eg:Nash} show that, even in the simple case of functions over $R$, we have the following strict inclusions:
		\begin{enumerate}
			\item $\ZC^\infty(R)\subsetneq\dots \subsetneq \ZC^{k+1}(R)\subsetneq\ZC^{k}(R)\subsetneq\dots\subsetneq\ZCc(R)$;
			\item $\mathcal{R}(R)\subsetneq \ZC^\infty(R)\subsetneq \mathcal{N}(R)$.
		\end{enumerate}
		
		\begin{rem}
			{Denote by $\mathcal{R}^k(\mathbb{R}^n)$ be the ring of $k$-regulous functions in the sense of \cite{FHMM16}. As the indeterminacy locus of the rational representation of a $k$-regulous function $f\in\mathcal{R}^k(\mathbb{R}^n)$ for $k\in\mathbb{N}$ is an algebraic subset of $R^n$ of codimension at least two by \cite[Prop.\,3,5]{FHMM16}, we have that $\mathcal{R}^k(\mathbb{R})=\mathcal{R}(\mathbb R)$ for every $k\in\mathbb{N}$. Moreover, $\mathcal{R}^\infty(\mathbb{R}^n)=\mathcal{R}(\mathbb R^n)$ for every $n\in\mathbb{N}$ by \cite[Thm.\,3.3]{FHMM16}. This shows that the ring $\ZC^k(\mathbb R)$ extends non-trivially even the ring $\mathcal{R}^k(\mathbb{R})$ of $k$-regulous functions for every $k\in\mathbb{N}\cup\{\infty\}$. We will see in Lemma \ref{lem:regulous} below that this is actually a general behavior.}
		\end{rem}

		\subsubsection{The $\ZC^k$-Łojasiewicz inequality} The next result is a $\ZC^k$-version of the Łojasiewicz inequality when $k\in\mathbb{N}$ which will turn out to be a very useful tool in Section \ref{subsec:ZCk(Rn)}. Let us fix the following notation: let $U\subset R^n$ be a $\ZC$-open set and let $f\in\ZC^k(U)$, then denote by $\mathcal{Z}(f)\coloneq \{x\in U : f(x)=0\}$ the \emph{zero locus} of $f$ in $U$.
		
		\begin{prop}\label{prop:LojP}
			Let $U\subset R^n$ be a $\ZC$-open subset, let $f\in\ZC^k(U)$ and $g\in\ZC^k(U\setminus \mathcal{Z}(f))$. Then there exists $N\in\mathbb N\setminus\{0\}$ such that $h\in\ZC^k(U)$ satisfying $h=f^Ng$ on $U\setminus \mathcal{Z}(f)$ and $h=0$ on $\mathcal{Z}(f)$.
		\end{prop}
		\begin{proof}
			As $U\subset R^n$ is a $\ZC$-open subset of $R^n$ it is in particular a locally closed semialgebraic set, therefore by \cite[Proposition 2.6.4]{BCR}, there exists $M\in\mathbb N\setminus\{0\}$ such that for every $\alpha\in\mathbb{N}^n$ with $|\alpha|\leq k$ there is a semialgebraic and continuous function $h_\alpha:U\to R$ such that $h_\alpha|_{U\setminus \mathcal{Z}(f)}=f^M \frac{\partial^{|\alpha|} g}{\partial x^\alpha}$ and $h_\alpha|_{\mathcal{Z}(f)}\equiv 0$. Fix $N\coloneq M+k$. Therefore, the function $h:U\to R$ defined as $h|_{U\setminus \mathcal{Z}(f)}=f^N g$  and $h|_{\mathcal{Z}(f)}\equiv 0$ is of class $\mathcal{C}^k$. Note that $\Gamma_h=\left(\Gamma_{f^Ng|_{U\setminus \mathcal{Z}(f)}}\right)\sqcup\left(\mathcal{Z}(f)\times\{0\}\right)$ is a $\ZC$-subset of $R^{n+1}$ as the union of two $\ZC$-sets. Hence, $h\in\ZC^k(U)$, as desired.
		\end{proof}

		\begin{rem}\label{rem:Loj-infty}
			{Observe that Proposition \ref{prop:LojP} does not hold for $k=\infty$. Consider $f\in\ZC^\infty(R)$ defined as $f(x)=x$ and $g\in\ZC^\infty(R\setminus\{0\})$ defined as $g(x)=\sqrt[3]{x}$, then $h_N\coloneq f^N g=x^{N}\sqrt[3]{x}\in\ZC^N(R)\setminus\ZC^{N+1}(R)$ for every $N\in\mathbb{N}$.}
		\end{rem}
		
		As a direct consequence of Proposition \ref{prop:LojP}, we get the following version of the Łojasiewicz inequality, which should be considered as a weak version of the Nullstellensatz we will prove in Theorem \ref{thm:Nullstellensatz}.
		
		\begin{cor}[Łojasiewicz inequality]\label{cor:Loj-ineq}
			Let $U\subset R^n$ be a $\ZC$-open set and let $f,g\in\ZC^k(U)$ for $k\in\mathbb{N}$. If $\mathcal{Z}(g)\subset \mathcal{Z}(f)$ then there exist $N\in\mathbb N\setminus\{0\}$ and $h\in\ZC^k(U)$ such that $f^N=hg$.
		\end{cor}
		\begin{proof}
			An application of Proposition \ref{prop:LojP} with the following substitutions: ``$f$"$\coloneq f$ and ``$g$"$\coloneq \frac1g\in\ZC^k(U\setminus \mathcal{Z}(f))$ gives $h\in\ZC^k(U)$ such that $h=\frac{f^N}{g}$ on $U\setminus\mathcal{Z}(f)$ and $h=0$ on $\mathcal{Z}(f)$. Hence $f^N=hg$ on $U$, as desired.
		\end{proof}

	In the rest of the paper we will mainly restrict to the case $k\in \mathbb{N}$, as the properties of the ring of $\ZC^k$-functions for $k\in\mathbb{N}$ and $k=\infty$ are very different, since the $\ZC^k$-Łojasiewicz inequality holds only when $k\in\mathbb{N}$, as mentioned in Remark \ref {rem:Loj-infty}. The properties of the ring $\ZC^\infty(U)$, with $U$ a $\ZC$-open subset of $R^n$ are the topic of an ongoing work.
	
	\subsubsection{Zero loci of $\ZC^k$-functions}\label{subsec:zeros}
	
	Next we show that the $\ZC$-closed subsets of $R^n$ actually coincide with the zero loci of $\ZC^k$-functions for $k\in\mathbb{N}$.
	
	\begin{prop}\label{prop:zerolocus}
		Let $X\subset R^n$ be a $\ZC$-closed set and $k\in\mathbb{N}$. Then there exists $f\in\ZC^k(R^n)$ such that $X=\mathcal{Z}(f)$. 
	\end{prop}
	\begin{proof}
		Let us prove the result by induction on $d\coloneq \dim X\geq 0$. If $d=0$, $X$ consists just of a finite number of points and the statement is trivial. Assume then that $d>0$ and the statement holds for every $\ZC$-closed subset $Y$ of ${R}^n$ of dimension $<d$.\\
		
		By Remark \ref{rem:irrcomps}, $X$ decomposes uniquely into a finite number of $\ZC$-closed irreducible components $X_1,\dots,X_s$. Thus, up to considering each $\ZC$-closed irreducible component of $X$ separately, we may suppose $X$ is an irreducible $\ZC$-closed set. Let $V\coloneq \zarcl{X}$ be the Zariski closure of $X$ in ${R}^n$ and let $f\in  R[x]$ such that $V=\mathcal{Z}(f)$. As $V$ is algebraic, it is in particular a $\ZC$-closed subset of ${R}^n$, hence it admits a unique decomposition $\mathcal{V}\coloneq \{V_1,\dots,V_t\}$ into $\ZC$-irreducible components subsets of ${R}^n$ by Remark \ref{rem:irrcomps}. Let $Z\coloneq \bigcup_{W\in\mathcal{V}, \dim W<d} W$. As $\dim Z<d$, let $g\in\ZC^k( R^n)$ such that $Z=\mathcal{Z}(g)$. Moreover, as $X$ is a $\ZC$-set, it is a finite union of Zariski locally closed subsets of $ R^n$, so in particular $\dim (V\setminus X)$<d. Thus, $V=X\cup Z$.\\
		
		Let $U\coloneq {R}^n\setminus X$. Observe that $\mathcal{Z}(f)\cap U=\mathcal{Z}(g)\cap U$ and $\frac{1}{f^2}\in\ZC^k(U\setminus \mathcal{Z}(g))$, thus an application of Proposition \ref{prop:LojP} to ``$U$"$\coloneq  U$, ``$f$"$\coloneq g$ and ``$g$"$\coloneq\frac{1}{f^2}$ gives $N\in \mathbb{N}\setminus\{0\}$ and $h\in\ZC^k(U)$ such that $h=\frac{g^N}{f^2}$ on $U\setminus \mathcal{Z}(g)=R^n\setminus(X\cup Z)=U\setminus V$ and $h\equiv 0$ on $U\cap \mathcal{Z}(g)=Z\setminus X$. Observe that the function $\frac{f^2(x)}{f^2(x)+g^{2N}(x)}\in\ZC^k( R^n\setminus Z)$, by Proposition \ref{prop:composition} and Corollary \ref{cor:ZCk-ring}, and it extends of class $\mathcal{C}^k$ to $s: R^n\setminus (X\cap Z)\to {R}$ defined as
		\[
		s(x)\coloneq 
		\begin{cases}
			\frac{f^2(x)}{f^2(x)+g^{2N}(x)}\quad\quad &x\in  R^n\setminus Z,\\
			1  & x\in Z\setminus X.
		\end{cases}
		\]
		Observe that, as $R^n\setminus Z$ and $Z\setminus X$ are disjoint $\ZC$-sets such that $(R^n\setminus Z)\cup(Z\setminus X)= R^n\setminus (X\cap Z)$, $s: R^n\setminus (X\cap Z)\to  R$ is of class $\mathcal{C}^k$ by the choice of $N\in\mathbb{N}\setminus\{0\}$ and $\Gamma_{s}=\Gamma_{s|_{R^n\setminus Z}}\sqcup\left((Z\setminus X)\times\{1\}\right)\subset R^{n+1}$, therefore $s\in\ZC^k(R^n\setminus (X\cap Z))$. Moreover, $\mathcal{Z}(s)=(R^n\setminus Z)\cap V=X\setminus Z$.\\
		
		By definition of $Z\coloneq \bigcup_{W\in\mathcal{V}, \dim W<d} W$, we have that $\dim (X\cap Z)<\dim X\eqcolon d$, then, by inductive assumption, there is $\ell\in\ZC^k({R}^n)$ such that $\mathcal{Z}(\ell)=X\cap Z$. Therefore, an application of Proposition \ref{prop:LojP} to ``$U$"$\coloneq  {R}^n$, ``$f$"$\coloneq\ell$ and ``$g$"$\coloneq s\in\ZC^k(R^n\setminus\mathcal{Z}(\ell))$ gives $M\in \mathbb{N}\setminus\{0\}$ and $t\in\ZC^k({R}^n)$ such that $t|_{{R}^n\setminus(X\cap Z)}=(\ell|_{{R}^n\setminus(X\cap Z)})^M s$ and $t|_{X\cap Z}\equiv 0$. Hence, $t\in\ZC^k( R^n)$ satisfies $\mathcal{Z}(t)=\mathcal{Z}(s)\cup (X\cap Z)=(X\setminus Z)\cup(X\cap Z)=X$, as desired.
	\end{proof}

\subsection{The ring of $\ZC^k$-functions}\label{subsec:ZCk(Rn)}

	\subsubsection{The Nullstellensatz and radical ideals of $\ZC^k(R^n)$}\label{subsec:Nullstellensatz}
	
	Let us introduce the following notation, similar to the usual one for zero sets and vanishing ideals in algebraic geometry. Let $k\in\mathbb{N}\cup\{\infty\}$. Let $S\subset \ZC^k(R^n)$ and $X\subset R^n$, we denote by
	\begin{align*}
		\mathcal{Z}(S)&\coloneq\{x\in R^n\,:\,f(x)=0\text{ for every $f\in S$}\},\\
		\I_{\ZC^k(R^n)}(X)&\coloneq\{f\in\ZC^k(R^n)\,:\,f(x)=0\text{ for every $x\in X$}\}
	\end{align*}
	the \emph{vanishing set} $\mathcal{Z}(S)$ of $S$ in $R^n$ and the \emph{vanishing ideal} $\I_{\ZC^k(R^n)}(X)$ of $X$ in $\ZC^k(R^n)$. When $S=\{f_1,\dots,f_s\}$ is finite, we shorten the notation as $\mathcal{Z}(f_1,\dots,f_s)\coloneq \mathcal{Z}(\{f_1,\dots,f_s\})$. Moreover, $\mathcal{Z}(\I_{\ZC^k(R^n)}(\mathcal{Z}(S)))=\mathcal{Z}(S)$ for every $S\subset \ZC^k(R^n)$.\\
	
	We stress that, although the $\ZC$-topology of $R^n$ is Noetherian, the ring $\ZC^k(R^n)$ is not Noetherian for any $k\in\mathbb N$ as the following example shows.
	
	\begin{eg}\label{eg:notnoeth}
		Let $k\in\mathbb N$. Consider the ideal $I\coloneq (f_r)_{r\in\mathbb{N}\setminus\{0\}}$ of $\ZC^k(R^n)$ where $f_r:R^n\to R$ is defined as
		\[
		f_r(x)\coloneq x_1^k\sqrt[2r+1]{x_1},
		\]
		for any $r\in\mathbb N\setminus\{0\}$. Let us show that $I$ is not finitely generated. Let us show that $f_{m+1}\notin (f_1,\dots,f_m)$ for any $m\in\mathbb{N}\setminus\{0\}$. Assume  by contradiction that there are $h_1,\dots,h_m\in\ZC^k(R^n)$ such that $f_{m+1}=\sum_{r=1}^{m} h_r f_r$. Observe that, by definition of $f_r$ for $r\in\mathbb{N}\setminus\{0\}$, this is equivalent to have that $\sqrt[2m+3]{x}=\sum_{r=1}^{m} h_r \sqrt[2r+1]{x}$. In particular, by continuity at $0$ of all the involved functions, this implies that $1=\lim_{x\to 0} \frac{\sum_{r=1}^{m} h_r(x) \sqrt[2r+1]{x}}{\sqrt[2m+3]{x}}=\sum_{r=1}^{m}\lim_{x\to 0} h_r(x) \sqrt[(2m+3)(2r+1)]{x^{2(m-r+1)}}$. Again, by continuity at $0$ (as $r\in\{1,\dots,m\}$) we deduce that $\lim_{x\to 0} h_r(x) \sqrt[(2m+3)(2r+1)]{x^{2(m-r+1)}}=0$ for every $r\in\{1,\dots,m\}$, getting to a contradiction.
	\end{eg} 
	
	In what follows, we replace the usual arguments relying on the finite generation of ideals with the Noetherianity of the $\ZC$-topology and the Łojasiewicz inequality.
	
	\begin{lemma}\label{lem:Iprincipal}
		Let $k\in\mathbb{N}\cup\{\infty\}$. Let $I\subset \ZC^k(R^n)$ be an ideal, then there exists $f\in I$ such that $\mathcal Z(f)=\mathcal Z(I)$.
	\end{lemma}
	\begin{proof}
		By Noetherianity of the $\ZC$-topology of $R^n$ there are $g_1,\ldots,g_r\in I$ such that $\mathcal Z(I)\coloneq\bigcap_{g\in I}\mathcal Z(g)=\bigcap_{i=1}^r\mathcal Z(g_i)$. Therefore, $f\coloneq\sum_{i=1}^rg_i^2\in I$ satisfies $\mathcal{Z}(f)=\mathcal Z(I)$, as desired.
	\end{proof}

\begin{rem}\label{rem:ZI-closed}
Let $k\in\mathbb{N}$. By Lemma \ref{lem:Iprincipal} and Remark \ref{rmk:zeros}, $\mathcal{Z}(I)$ is a $\ZC$-closed subset of $R^n$ for every ideal $I\subset\ZC^k(R^n)$ (this part also holds for $k=\infty$). Conversely, by Proposition \ref{prop:zerolocus}, every $\ZC$-closed subset $X$ of $R^n$ satisfies $X=\mathcal{Z}(f)$ for some $f\in\ZC^k(R^n)$. Hence, when $k\in\mathbb{N}$, the sets $\mathcal{Z}(I)$, where $I$ are ideals of $\ZC^k(R^n)$, are exactly the $\ZC$-closed subsets of $R^n$.
\end{rem}

	\begin{prop}[Weak Nullstellensatz]\label{prop:w-Nullstellensatz}
		Let $k\in\mathbb{N}\cup\{\infty\}$. Let $I\subset\ZC^k(R^n)$ be an ideal. If $\mathcal{Z}(I)=\varnothing$, then $I=\ZC^k(R^n)$.
	\end{prop}
	\begin{proof}
	    By Lemma \ref{lem:Iprincipal}, there exists $f\in I$ such that $\mathcal{Z}(f)=\mathcal{Z}(I)=\varnothing$. By Proposition \ref{prop:composition}, $\frac1f\in\ZC^k(R^n)$. Hence $1=\frac1f\cdot f\in I$.
	\end{proof}

	\begin{thm}[Nullstellensatz]\label{thm:Nullstellensatz}
		Let $k\in\mathbb{N}$. Let $I\subset\ZC^k(R^n)$ be an ideal. Then
		\[
		\I_{\ZC^k(R^n)}(\mathcal{Z}(I))=\sqrt{I}.
		\]
	\end{thm}
	\begin{proof}
		By the above lemma, there exists $f\in I$ such that $\mathcal Z(f)=\mathcal Z(I)$. Let $g\in\I_{\ZC^k(R^n)}(\mathcal{Z}(I))$ then $\mathcal Z(f)=\mathcal Z(I)\subset\mathcal Z(g)$. By the Łojasiewicz inequality (see Corollary \ref{cor:Loj-ineq}), there is $N\in\mathbb N\setminus\{0\}$ and $h\in\ZC^k(R^n)$ such that $g^N=hf\in I$. This shows that $g\in\sqrt{I}$, as desired. The other inclusion is obvious.
	\end{proof}

	Using Remark \ref{rem:ZI-closed}, we can restate the Nullstellensatz as a correspondence between radical ideals of $\ZC^k(R^n)$ and $\ZC$-closed subsets of $R^n$ as follows.
	\begin{cor}\label{cor:rad-closed}
		Let $k\in\mathbb N$. The maps $I\mapsto\mathcal{Z}(I)$ and $X\mapsto\I_{\ZC^k(R^n)}(X)$ define inclusion reversing bijections between the set of radical ideals of $\ZC^k(R^n)$ and the set of $\ZC$-closed subsets of $R^n$ which are inverse of each other.
		In particular, if $I,J\subset\ZC^k(R^n)$ are ideals then $\sqrt{I}\subset\sqrt{J}$ if and only if $\mathcal{Z}(J)\subset\mathcal{Z}(I)$.
	\end{cor}

	The next result shows that the ring $\ZC^k(R^n)$ is \emph{radically principal} for every $k\in\mathbb N$, that is, for every ideal $I\subset\ZC^k(R^n)$ there exists $f\in I$ such that $\sqrt{I}=\sqrt{(f)}$.
	
	\begin{thm}\label{thm:radprinc}
		The ring $\ZC^k(R^n)$ is radically principal for every $k\in\mathbb N$.
	\end{thm}
	\begin{proof}
		Let $I\subset\ZC^k(R^n)$ be an ideal. By Lemma \ref{lem:Iprincipal}, there exists $f\in I$ such that $\mathcal Z(f)=\mathcal Z(I)$. By the Nullstellensatz (see Theorem \ref{thm:Nullstellensatz}), we have that $\sqrt{I}=\I_{\ZC^k(R^n)}(\mathcal Z(I))=\I_{\ZC^k(R^n)}(\mathcal Z(f))=\sqrt{(f)}$, as desired.
	\end{proof}
	
	The element $f\in I$ provided by the above Theorem \ref{thm:radprinc} is not unique, the next result characterizes the suitable $f$.
	
	\begin{cor}\label{cor:radgeneration&zeros}
		Let $k\in\mathbb N$. Let $I\subset\ZC^k(R^n)$ be an ideal and let $f\in I$. Then, $\sqrt{(f)}=\sqrt{I}$ if and only if $\mathcal{Z}(f)=\mathcal{Z}(I)$.
	\end{cor}
	\begin{proof}
		An application of Corollary \ref{cor:rad-closed} to the ideals $(f)$ and $I$ gives that $\sqrt{(f)}=\sqrt{I}$ if and only if $\mathcal{Z}((f))=\mathcal{Z}(I)$. Then, it suffices to observe that $\mathcal{Z}((f))=\mathcal{Z}(f)$.
	\end{proof}
	
	Recall that an ideal $I$ of a commutative ring $A$ is \emph{real} if whenever $f_1^2+\cdots+f_m^2\in I$ with $f_1,\ldots,f_m\in A$, then $f_i\in I$ for every $i\in\{1,\ldots,m\}$.
	
	\begin{cor}
		Let $I\subset\ZC^k(R^n)$ be an ideal for $k\in\mathbb N$. Then,  $I$ is radical if and only if $I$ is real.
	\end{cor}
	\begin{proof}
		A real ideal is always radical. 
		
		Let $I\subset\ZC^k(R^n)$ be a radical ideal. Let $f_1^2+\cdots+f_m^2\in I$ with $f_1,\ldots,f_m\in\ZC^k(R^n)$. We derive from the Nullstellensatz that $f_1^2+\cdots+f_m^2\in I=\I_{\ZC^k(R^n)}(\mathcal Z(I))$. Thus $f_1^2+\cdots+f_m^2$ vanishes on $\mathcal Z(I)$, hence each $f_i$ vanishes on $\mathcal Z(I)$, i.e. $f_i\in\I_{\ZC^k(R^n)}(\mathcal Z(I))=I$.
	\end{proof}
	
	For $k\in\mathbb N\cup\{\infty\}$ and $a=(a_1,\dots,a_n)\in R^n$, we denote by
	\begin{align*}
		\mathfrak{m}_a\coloneq\{f\in\ZC^k(R^n) : f(a)=0\}
	\end{align*}
	the maximal ideal of $\ZC^k$-functions vanishing at $a\in R^n$. As $\mathcal{Z}\left(\mathfrak{m}_a\right)=\mathcal{Z}(x_1-a_1,\dots,x_n-a_n)=\mathcal{Z}\left(\sum_{i=1}^{n}(x_i-a_i)^2\right)$, an application of Corollary \ref{cor:radgeneration&zeros} gives that
	\[
	\mathfrak{m}_a=\sqrt{(x_1-a_1,\dots,x_n-a_n)\cdot \ZC^k(R^n)}=\sqrt{\left(\sum_{i=1}^{n}(x_i-a_i)^2\right)\cdot \ZC^k(R^n)}.
	\]
	
	Let us characterize maximal ideals $\mathfrak{m}$ in $\ZC^k(R^n)$, with $k\in \mathbb{N}\cup\{\infty\}$, as follows. 
	
	\begin{cor}\label{cor:max-ideals}
		Let $k\in\mathbb N\cup\{\infty\}$ and let $\mathfrak{m}\subset\ZC^k(R^n)$ be an ideal. Then, $\mathfrak{m}$ is maximal if and only if $\mathfrak{m}=\mathfrak{m}_a$ for some $a\in R^n$.
	\end{cor}
	\begin{proof}
		It is clear that $\mathfrak{m}_a$ is maximal for every $a\in R^n$. Conversely, let $\mathfrak{m}$ be a maximal ideal of $\ZC^k(R^n)$. As $\mathfrak{m}\neq\ZC^k(R^n)$, we deduce from the weak Nullstellensatz (see Proposition \ref{prop:w-Nullstellensatz}) that $\mathcal{Z}(\mathfrak{m})\neq\varnothing$. Let $a\in\mathcal{Z}(\mathfrak{m})$. Then $\mathfrak{m}\subset\mathfrak{m}_a$, and since $\mathfrak{m}$ is maximal and $\mathfrak{m}_a\neq\ZC^k(R^n)$, we actually have $\mathfrak{m}=\mathfrak{m}_a$, as desired.
	\end{proof}

	In what follows, we denote by $\|x\|\coloneq\sqrt{\sum_{i=1}^nx_i^2}$ the Euclidean norm of $x\in R^n$.
	\begin{cor}
		Let $k\in\mathbb N$, $f\in\ZC^k(R^n)$ and $a\in\mathcal{Z}(f)$. Then there is $N\in \mathbb{N}\setminus\{0\}$ such that $f^N$ is differentiable at $a$ and $D_a(f^N)=0$.
	\end{cor}
	\begin{proof}
		Consider the polynomial $q(x)\coloneq\|x-a\|^2=\sum_{i=1}^n(x_i-a_i)^2$. As $\mathcal{Z}(q)=\{a\}\subset\mathcal{Z}(f)$, we have that $\frac1q\in\ZC^k(R^n\setminus\mathcal{Z}(f))$. An application of Proposition \ref{prop:LojP} with ``$U$''$\coloneq R^n$, ``$f$''$\coloneq f$ and ``$g$''$\coloneq\frac1q$ gives $N\in\mathbb N\setminus\{0\}$ and $h\in\ZC^k(R^n)$ such that $h=\frac{f^N}{q}$ on $R^n\setminus\mathcal{Z}(f)$ and $h\equiv0$ on $\mathcal{Z}(f)$. Therefore $f^N=hq$ on $R^n$, and $h(a)=0$ as $a\in\mathcal{Z}(f)$.
		
		Thus $|f^N(x)|=|h(x)|\cdot\|x-a\|^2=o\left(\|x-a\|^2\right)$ as $x\to a$. Since $f^N(a)=0$, this implies that $f^N$ is differentiable at $a$ with $D_a(f^N)=0$, as desired.
	\end{proof}
	
	\subsubsection{The Nullstellensatz on open sets}
	
	We obtain as well the following characterization of $\ZC^k$ functions on $\ZC$-open sets.
	
	\begin{prop}\label{prop:localization}
		Let $k\in\mathbb N$. Let $U\subset R^n$ be a $\ZC$-open set and let $f\in\ZC^k(R^n)$ such that $\mathcal{Z}(f)=R^n\setminus U$. Denote by $\rho:\ZC^k(R^n)\to \ZC^k(U)$ the restriction morphism $\rho(h)\coloneq h|_{U}$ for every $h\in\ZC^k(R^n)$. Then, the induced morphism:
		\[
		\rho_f:\ZC^k(R^n)_f\to\ZC^k(U),\quad \frac{h}{f^N}\mapsto \frac{h|_U}{f^N|_{U}}
		\]
		is a canonical isomorphism, that is, whenever ${g}\in\ZC^k(R^n)$ satisfies $\mathcal{Z}({g})=R^n\setminus U$ there is a unique isomorphism $\varphi:\ZC^k(R^n)_{f}\to\ZC^k(R^n)_g$ of $\ZC^k(R^n)$-algebras such that $\rho_{f}=\rho_g\circ\varphi$.
	\end{prop}

	\begin{proof}
		Let $g\in \ZC^k(U)$. By Proposition \ref{prop:LojP} applied with the following substitutions: ``$U$"$\coloneq  R^n$, ``$f$"$\coloneq f$ and ``$g$"$\coloneq g\in\ZC^k(R^n\setminus\mathcal{Z}(f))$ gives $N\in \mathbb{N}\setminus\{0\}$ and $h\in\ZC^k(R^n)$ such that $h=f^N g$ on $U={R}^n\setminus \mathcal{Z}(f)$ and $h\equiv 0$ on $\mathcal{Z}(f)$. Therefore, $\rho_f\left(\frac{h}{f^N}\right)=g$ showing that $\rho_f$ is surjective. Let us verify that $\rho_f$ is injective. Let $h\in\ZC^k(R^n)$ and $N\in\mathbb N$ be such that $\rho_f\left(\frac{h}{f^N}\right)=0$, i.e. $h_{|U}=0$. Then $fh=0$ on $R^n$ (since $h$ vanishes on $U$ and $f$ vanishes on $R^n\setminus U$). Hence $\frac{h}{f^N}=0$ in $\ZC^k(R^n)_f$.

Let $g\in\ZC^k(R^n)$ be such that $\mathcal{Z}(g)=R^n\setminus U$. By the above, $\rho_g:\ZC^k(R^n)_g\to\ZC^k(U)$ is an isomorphism so that $\varphi\coloneq\rho_g^{-1}\circ\rho_f$ is an isomorphism of $\ZC^k(R^n)$-algebras satisfying $\rho_f=\rho_g\circ\varphi$. It is unique by the universal property of the localization.
	\end{proof}
	
	\begin{rem}
		Observe that Proposition \ref{prop:localization} does not hold in the classical case where $U\subset R^2$ is a Zariski-open set and the considered functions are regular instead of $\ZC^k$. Indeed, consider $f\in R[x,y]$ as $f(x,y)=x$, the open set $U\coloneq R^2\setminus\mathcal{Z}(f)$ and $g\in \mathcal{R}(U)$ defined as $g(x,y)\coloneq \frac{1}{x^2+y^2}$. Denote by $\rho:\mathcal{R}(R^2)\to \mathcal{R}(U)$ defined as $\rho(h)=h|_U$ for every $h\in\mathcal{R}(R^2)$ the restriction morphism and by $\rho_f:\mathcal{R}(R^2)_f\to \mathcal{R}(U)$ the induced morphism. Evidently, $\I(\mathcal{Z}(f))=(f)$ and $g=x^2+y^2\notin (x)=(f)$, therefore there do not exist $h\in\mathcal{R}(R^2)$ and $N\in\mathbb{N}$ such that $\rho_f\left(\frac{h}{x^N}\right)=g$. This shows that $\rho_f$ is not surjective.
	\end{rem}

	Let us extend the notation of Subsection \ref{subsec:Nullstellensatz} over $\ZC$-open sets. Let $k\in\mathbb{N}\cup\{\infty\}$. Let $U\subset R^n$ be a $\ZC$-open set. Let $S\subset\ZC^k(U)$ and $X\subset U$. Denote by
	\begin{align*}
		\mathcal{Z}(S)&\coloneq \{x\in U\,:\,f(x)=0 \text{ for all $f\in S$}\}\\
		\I_{\ZC^k(U)}(X)&\coloneq \{f\in\ZC^k(U)\,:\,f(x)=0 \text{ for every $x\in X$}\}
	\end{align*}
	the \emph{vanishing set} $\mathcal{Z}(S)$ of $S$ in $U$ and the \emph{vanishing ideal} $\I_{\ZC^k(U)}(X)$ of $X$ in $\ZC^k(U)$. When $S=\{f_1,\dots,f_s\}$ is finite, we shorten the notation as $\mathcal{Z}(f_1,\dots,f_s)\coloneq \mathcal{Z}(\{f_1,\dots,f_s\})$. Moreover, $\mathcal{Z}(\I_{\ZC^k(U)}(\mathcal{Z}(S)))=\mathcal{Z}(S)$ for every $S\subset \ZC^k(U)$.  Observe that, by Remark \ref{rmk:zeros}, $\mathcal{Z}(S)=\bigcap_{f\in S}\mathcal{Z}(f)$ is closed in $U$ for the induced $\ZC$-topology, so that $\mathcal Z(S)$ is $\ZC$-locally closed. \\
	
	Let $U\subset R^n$ be a $\ZC$-open set. Here we provide a version of the Nullstellensatz for the ring $\ZC^k(U)$, when $k\in\mathbb{N}$.
	
	\begin{cor}[Nullstellensatz for $\ZC^k(U)$]\label{cor:Nullstellensatz-U}
		Let $k\in\mathbb{N}$. Let $U\subset R^n$ be a $\ZC$-open set and let $I\subset\ZC^k(U)$ be an ideal. Then
		\[
		\I_{\ZC^k(U)}(\mathcal{Z}(I))=\sqrt{I}.
		\]
	\end{cor}
	
	\begin{proof}
		Let $f\in\ZC^k(R^n)$ such that $\mathcal{Z}(f)=R^n\setminus U$. By Proposition \ref{prop:localization},
		the morphism $\rho_f:\ZC^k(R^n)_f\to \ZC^k(U)$ defined as $\rho_f\left(\frac{h}{f^N}\right)=\frac{h|_U}{f^N|_{U}}$ is an isomorphism. Denote by $J\coloneq \rho_f^{-1}(I)$. As $J=(J\cap \ZC^k(R^n))\cdot \ZC^k(R^n)_f$, an application of the Nullstellensatz (see Theorem \ref{thm:Nullstellensatz}) to the ideal $J\cap \ZC^k(R^n)$ of $\ZC^k(R^n)$ gives
		\begin{align*}
			\sqrt{I}&=\rho_f\left(\sqrt{J}\right)=\rho_f\left(\sqrt{J\cap\ZC^k(R^n)}\cdot\ZC^k(R^n)_f\right)\\
			&=\rho_f\left(\I_{\ZC^k(R^n)}(\mathcal{Z}(J\cap\ZC^k(R^n)))\cdot\ZC^k(R^n)_f\right)=\I_{\ZC^k(U)}(\mathcal{Z}(I)),
		\end{align*}
		as desired.
	\end{proof}
	
	As a consequence, we get that $\ZC^k(U)$ is radically principal as well.
	
	\begin{cor}\label{cor:radprinc-U}
		The ring $\ZC^k(U)$ is radically principal for every $k\in\mathbb N$.
	\end{cor}
	\begin{proof}
		Let $I\subset\ZC^k(U)$ be an ideal. By Proposition \ref{prop:localization}, the morphism $\rho_f:\ZC^k(R^n)_f\to \ZC^k(U)$ defined as $\rho_f\left(\frac{h}{f^N}\right)=\frac{h|_U}{f^N|_U}$ is an isomorphism. Denote by $J\coloneq \rho_f^{-1}(I)$. As $J=(J\cap \ZC^k(R^n))\cdot \ZC^k(R^n)_f$, an application of Theorem \ref{thm:radprinc} to the ideal $J\cap \ZC^k(R^n)$ of $\ZC^k(R^n)$ gives $h\in\ZC^k(R^n)$ such that $\sqrt{(h)}=\sqrt{J\cap \ZC^k(R^n)}$. Therefore,
		\[
		\sqrt{I}=\rho_f(J)=\rho_f\left(\sqrt{J\cap\ZC^k(R^n)}\cdot\ZC^k(R^n)_f\right)=\rho_f\left(\sqrt{(h)}\cdot\ZC^k(R^n)_f\right)=\sqrt{(h|_{U})},
		\]
		as desired.
	\end{proof}

		\subsubsection{$\ZC^k$-functions on $\ZC$-locally closed sets \& extension}
		
        Let $X\subset R^n$. Denote by $\zcl{X}$ the closure of $X$ in the $\ZC$-topology of $R^n$. We say that a $\ZC$ subset $X$ of $R^n$ is \emph{$\ZC$-locally closed} if there exists a $\ZC$-open subset $U$ of $R^n$ such that $X=\zcl{X}\cap U$. Note that, in this case, $\zcl{X}\setminus X=\zcl{X}\setminus U$ is $\ZC$-closed (and does not depend on the choice of $U$). In particular, by Proposition \ref{prop:zerolocus}, for every $k\in\mathbb{N}$ there exists $f\in\ZC^k(R^n)$ such that $\mathcal{Z}(f)=\zcl{X}\setminus X$.
		
\begin{rem}
	A $\ZC$-subset of $R^n$ which is locally closed for the Euclidean topology need not be $\ZC$-locally closed. Consider for instance the $\ZC$-set
	\[
	X\coloneq \{(x,y,z)\in R^3\,:\, zx^2=y^2,\ x^2+y^2\neq 0\}\cup\{(0,0,-1)\},
	\]
	that is, the Whitney umbrella with the $z$-axis removed (see Remark \ref{rem:whitney}) to which we add the point $(0,0,-1)$. Besides, $\ecl{X}=X\cup \{(0,0,z)\in R^3\,:\,z\ge 0\}$, and hence $X=\ecl{X}\cap\left(R^3\setminus\{(0,0,z)\in R^3\,:\,z\ge 0\}\right)$ is Euclidean locally closed. However, $\zcl{X}=\{(x,y,z)\in R^3\,:\,zx^2=y^2\}$ is the Whitney umbrella by Proposition \ref{prop:closure} and $\zcl{X}\setminus X=\{(0,0,z)\in\mathbb R^3\,:\,z\neq-1\}$ is not even Euclidean closed. Therefore, $X$ is not $\ZC$-locally closed.
\end{rem}
		
		Let us define the notion of $\ZC^k$-function on a $\ZC$-locally closed subset $X$ of $R^n$.
		
		\begin{defn}\label{def:ZC^k(X)}
			Let $k\in\mathbb{N}\cup\{\infty\}$. Let $X\subset R^n$ be a $\ZC$-locally closed set and let $f:X\to R$ be a function. We say that $f$ is a \emph{$\ZC^k$-function} if for every $x\in X$ there is a $\ZC$-open neighborhood $U_x$ of $x$ in $R^n$ and $F\in\ZC^k(U_x)$ such that $f|_{X\cap U_x}=F|_{X\cap U_x}$. We will denote by $\ZC^k(X)$ the set of $\ZC^k$-functions on $X$.
		\end{defn}
		
		\begin{rem}\label{rem:ZC-continuity}
			Let $X\subset R^n$ be a $\ZC$-locally closed set and let $f\in\ZC^k(X)$. Then:
			\begin{enumerate}[label=(\roman*),nosep]
				\item\label{rem:ZC-continuity-1} $f\in\ZC(X)$, indeed, by Noetherianity of the $\ZC$-topology of $R^n$ there are $\ZC$-open subsets $U_1,\dots,U_\ell$ of $R^n$ and $F_i\in\ZC^k(U_i)$ for every $i\in\{1,\dots,\ell\}$ such that $f|_{X\cap U_i}=F_i|_{X\cap U_i}$ for every $i\in\{1,\dots,\ell\}$ and $X\subset\bigcup_{i=1}^\ell U_i$. Therefore, $\Gamma_f=\bigcup_{i=1}^\ell \Gamma_{F_i}\cap (X\times R)$ is a $\ZC$-subset of $R^{n+1}$.
				\item\label{rem:ZC-continuity-2} As every $f\in\ZC^k(X)$ is continuous with respect to the Euclidean topology, then $f$ is continuous with respect to the $\ZC$-topology of $X$ by \ref{rem:ZC-continuity-1} and Proposition \ref{prop:ZCf}\ref{ZCf-2}.
				\item Let $U\subset R^n$ be a $\ZC$-open set, therefore $U$ is in particular $\ZC$-locally closed. In this case the notion of $\ZC^k$-function introduced in Definition \ref{def:ZC^k(X)} coincides with the one already introduced in Definition \ref{def:ZCk}. Indeed, suppose $f:U\to R$ satisfies Definition \ref{def:ZC^k(X)}. Then, $f\in\ZC(U)$ by \ref{rem:ZC-continuity-1} and $f$ is trivially of class $\mathcal{C}^k$ on the whole $U$.
			\end{enumerate}
		\end{rem}
		
		Let $X\subset R^n$ be a $\ZC$-locally closed set and let $f\in\ZCc(X)$. Observe that, $f\in\ZC(X)$ and $f$ is continuous by Remark \ref{rem:ZC-continuity}. The following open question then arises:
		
		\begin{question}
			Let $X\subset R^n$ be a $\ZC$-locally closed set. Let $f\in\ZC(X)$ be a continuous function, does it hold that $f\in\ZCc(X)$? In other words, for every $x\in X$, are there a $\ZC$-open neighborhood $U_x$ of $x$ in $R^n$ and a function $F\in\ZCc(U_x)$ such that $f|_{X\cap U_x}=F|_{X\cap U_x}$?
		\end{question}
		\noindent In the $k$-regulous setting, the analogous question has the following positive answers over~$\mathbb{R}$:
		\begin{itemize}[nosep]
			\item When $k=0$, continuous hereditarily rational functions on a (possibly singular) algebraic set $X\subset\mathbb R^n$ extend to regulous functions on $\mathbb R^n$, see \cite{KN15}.
			\item When $k>0$, $\mathcal{C}^k$ rational functions on a nonsingular algebraic set $X\subset\mathbb R^n$ in the sense of \cite{K09} extend to $k$-regulous functions on $\mathbb R^n$, see \cite{Ban24}. \\
		\end{itemize}
		
		In the next lemma we observe the relation between $\ZC^k$-functions and $k$-regulous functions in the sense of \cite{FHMM16}. Let $k\in\mathbb{N}\cup\{\infty\}$. A function $F:R^n\to R$ of class $\mathcal{C}^k$ is \emph{$k$-regulous} if there exist a Zariski open dense subset $U$ of $R^n$ and polynomials $p,q\in R[x_1,\dots,x_n]$, with $\mathcal{Z}(q)\cap U=\varnothing$, such that $F(x)=\frac{p(x)}{q(x)}$ for every $x\in U$. Let $X\subset R^n$ be $\ZC$-closed set. A function $f:X\to R$ is \emph{$k$-regulous} if there exists a $k$-regulous function $F:R^n\to R$ such that $f=F|_X$. Denote by $\mathcal{R}^k(X)$ the ring of $k$-regulous functions on $X$. 
		
		\begin{lemma}\label{lem:regulous}
			Let $X\subset \mathbb{R}^n$ be a $\ZC$-closed set and $f\in \mathcal{R}^k(X)$ be a $k$-regulous function with $k\in\mathbb{N}\cup\{\infty\}$, then $f\in\ZC^k(X)$.
		\end{lemma}
		
		\begin{proof}
			Recall that $\mathcal{R}^\infty(X)=\mathcal{R}(X)$, where $\mathcal{R}(X)$ denotes the ring of regular functions on $X$, therefore $\mathcal{R}^\infty(X)\subset\ZC^\infty(X)$. Let $k\in\mathbb{N}$. Let $F\in\mathcal{R}^k(\mathbb{R}^n)$ such that $f=F|_{X}$. By \cite[Rmk.\,5.39]{FHMM16}, $F$ is hereditarily rational in the sense of \cite{KN15}, that is, there are $0<\ell\leq n$ and a sequence of algebraic sets $X_0=X\supset X_1\supset \dots\supset X_\ell=\varnothing$ such that $f|_{X_i\setminus X_{i+1}}$ is a regular function for every $i\in\{0,\dots,\ell-1\}$. As each set $X_i\setminus X_{i+1}$ is a $\ZC$-subset of $\mathbb{R}^n$ and $f|_{X_i\setminus X_{i+1}}$ is regular, we have that $\Gamma_f=\bigcup_{i=0}^{\ell-1} \Gamma_{f|_{X_i\setminus X_{i+1}}}$ is a $\ZC$-subset of $\mathbb{R}^n$ as it is a finite union of $\ZC$-sets. Therefore, $F\in\ZC^k(\mathbb{R}^n)$ and hence $f=F|_X\in\ZC^k(X)$, as desired.
		\end{proof}
		
		Let $X\subset R^n$ be a $\ZC$-locally closed set. Let us show that $\ZC^k(X)$ has a ring structure.
		
		\begin{prop}\label{prop:ZC^k(X)-ring}
			Let $k\in\mathbb N\cup\{\infty\}$. Let $X\subset R^n$ be a $\ZC$-locally closed set. The set $\ZC^k(X)$ of $\ZC^k$-functions on $X$ endowed with the pointwise addition and multiplication is a ring. 
		\end{prop}
		\begin{proof}
			Let $f,g\in\ZC^k(X)$ and let $x\in X$. According to Definition \ref{def:ZC^k(X)}, there are $\ZC$-open neighborhoods $U_x,V_x$ of $x$ in $R^n$, $F\in\ZC^k(U_x)$ and $G\in\ZC^k(V_x)$ such that $f|_{X\cap U_x}=F|_{X\cap U_x}$ and $g|_{X\cap V_x}=G|_{X\cap V_x}$. Consider the $\ZC$-open neighborhood $U_x\cap V_x$ of $x$ in $R^n$ and the functions $F'\coloneq F|_{U_x\cap V_x}\in\ZC^k(U_x\cap V_x)$ and $G'\coloneq G|_{U_x\cap V_x}\in\ZC^k(U_x\cap V_x)$. By Corollary \ref{cor:ZCk-ring}, $F'+G',F'G'\in\ZC^k(U_x\cap V_x)$ and they satisfy $(F'+G')|_{U_x\cap V_x\cap X}=(f+g)|_{U_x\cap V_x\cap X}$ and $(F'G')|_{U_x\cap V_x\cap X}=(fg)|_{U_x\cap V_x\cap X}$. This shows that $f+g,fg\in\ZC^k(X)$.
		\end{proof}

		\begin{thm}[Extension of $\ZC^k$-functions]\label{thm:extension}
			Let $k\in\mathbb N$. Let $X\subset R^n$ be a $\ZC$-locally closed set and let $U\subset R^n$ be a $\ZC$-open set such that $X=U\cap\zcl{X}$. If $f\in\ZC^k(X)$, then there exists $F\in\ZC^k(U)$ such that $F|_{X}=f$.
		\end{thm}
		
		\begin{proof}
			Denote by $I$ the ideal of $\ZC^k(U)$ defined as:
			\[
			I\coloneq \{G\in\ZC^k(U)\,:\,  \text{there exists $H\in\ZC^k(U)$ such that $H|_{X}=f G|_{X}$}\}.
			\]
			It suffices to show that $1\in I$. According to the Nullstellensatz (see Corollary \ref{cor:Nullstellensatz-U}), it is enough to show that $\mathcal{Z}(I)=\varnothing$. Observe that choosing $g\in\ZC^k(R^n)$ such that $X=\mathcal{Z}(g)\cap U$ (see Proposition \ref{prop:zerolocus}), we can sett $G\coloneq g|_U$ and $H\coloneq G$ so that $\mathcal Z(I)\subset\mathcal Z(g)\cap U= X$. \\
			
			Let $x\in X$. By definition of $\ZC^k$ functions on $X$, there is a $\ZC$-open neighborhood $U_x$ of $x$ in $R^n$ and $F\in\ZC^k(U_x)$ such that $F|_{X\cap U_x}=f|_{X\cap U_x}$. Up to substituting $U_x$ by $U_x\cap U$ and $F$ by $F|_{U_x\cap U}$, we may suppose that $U_x\subset U$. 
			By Proposition \ref{prop:zerolocus}, there exists $G\in\ZC^k(R^n)$ such that $\mathcal{Z}(G)=R^n\setminus U_x$.
			An application of Proposition \ref{prop:LojP} with the following substitutions: ``$U$"$\coloneq U$, ``$f$"$\coloneq G|_U$ and ``$g$"$\coloneq F\in\ZC^k(U\setminus\mathcal{Z}(G))=\ZC^k(U_x)$ gives $N\in\mathbb N\setminus\{0\}$ and $H\in\ZC^k(U)$ such that $H=F G^N$ on $U_x$ and $H=0$ on $U\setminus U_x$.
			Therefore $H=fG^N$ on $X\cap U_x$, since $F=f$ on $X\cap U_x$, and $H=0=fG^N$ on $X\setminus U_x$, since $G$ vanishes on $R^n\setminus U_x$. Hence $H|_X=f\,G^N|_X$, that is, $G^N|_U\in I$. As $G(x)\neq0$, we get $x\notin\mathcal{Z}(I)$.\\
			
			This shows that $\mathcal{Z}(I)=\varnothing$, as desired.
		\end{proof}
		
		\begin{rem}
			Let $X\subset R^n$ be a $\ZC$-locally closed set and let $f\in\ZC^k(X)$. Then:
			\begin{enumerate}[label=(\roman*),nosep]
                \item As noticed at the beginning of this subsection, $\zcl{X}\setminus X$ is a $\ZC$-closed subset of $R^n$, therefore we can fix $U\coloneq R^n\setminus\left(\zcl{X}\setminus X\right)$ in the statement of Theorem \ref{thm:extension}. In particular, if $X\subset R^n$ is $\ZC$-closed we can fix $U\coloneq R^n$.
				\item If $X\subset R^n$ is not $\ZC$-closed, in general a function $f\in\ZC^k(X)$ cannot be extended to the whole $R^n$. Consider for instance $X=R\setminus\{0\}$ and $f:R\setminus\{0\}\to R$ defined by $f(x)=\frac1x$.\\
			\end{enumerate}
		\end{rem}
		
		In light of Proposition \ref{prop:localization}, we can then restate Theorem \ref{thm:extension} as follows.
		
\begin{cor}\label{cor:extension}
    Let $k\in\mathbb N$. Let $X\subset R^n$ be a $\ZC$-locally closed set. Let $f\in\ZC^k(R^n)$ be such that $\mathcal{Z}(f)\cap\zcl{X}=\zcl{X}\setminus X$. Denote by $I\coloneq\I_{\ZC^k(R^n)}(X)$, by $\bar{f}\in\ZC^k(R^n)/I$ the image of $f$ under the quotient map $\ZC^k(R^n)\to\ZC^k(R^n)/I$ and by $\nu:\ZC^k(R^n)/I\to \ZC^k(X)$ the morphism defined as $\nu(h+I)\coloneq h|_{X}$. Then the induced morphism
    \[
    \nu_f:\left(\ZC^k(R^n)/I\right)_{\overline{f}}\to \ZC^k(X),\quad \frac{h+I}{\bar{f}^N}\mapsto \frac{h|_{X}}{f^N|_{X}},
    \]
    is a canonical isomorphism in the following sense: let $g\in\ZC^k(R^n)$ be such that $\mathcal{Z}(g)\cap\zcl{X}=\zcl{X}\setminus X$, then there is a unique isomorphism $\varphi:\left(\ZC^k(R^n)/I\right)_{\bar{f}}\to \left(\ZC^k(R^n)/I\right)_{\bar{g}}$ of $\ZC^k(R^n)/I$-algebras such that $\nu_f=\nu_g\circ\varphi$.
\end{cor}
		
		\begin{proof}
			Denote by $U\coloneq R^n\setminus\mathcal{Z}(f)$ and denote by $J\coloneq \I_{\ZC^k(U)}(X)$. By Theorem \ref{thm:extension}, the restriction map $\psi:\ZC^k(U)/J\to\ZC^k(X)$ defined as $\psi(h+J)=h|_{X}$ for every $h\in\ZC^k(U)$ is an isomorphism. Denote by $\rho:\ZC^k(R^n)\to \ZC^k(U)$ the restriction map. By Proposition \ref{prop:localization}, the induced map $\rho_f:\ZC^k(R^n)_f\to \ZC^k(U)$ defined as $\rho_f\left(\frac{h}{f^N}\right)=\frac{h|_{U}}{f^N|_U}$ for every $h\in\ZC^k(R^n)$ is a (canonical) isomorphism. Therefore, the induced map $\widetilde{\rho_f}:\ZC^k(R^n)_f/I_f\to \ZC^k(U)/J$ is an isomorphism, where $I_f\coloneq I\cdot\ZC^k(R^n)_f$. Denote by $\phi:(\ZC^k(R^n)/I)_{\bar{f}}\to \ZC^k(R^n)_f/I_f$ the (canonical) isomorphism defined as $\phi\left(\frac{h+I}{\bar{f}^N}\right)\coloneq\left(\frac{h}{f^N}+I_f\right)$ for every $h\in\ZC^k(R^n)$. Then, $\nu_f=\psi\circ \widetilde{\rho_f}\circ \phi$ showing that $\nu_f:(\ZC^k(R^n)/I)_{\bar{f}}\to \ZC^k(X)$ is an isomorphism.
			Let $g\in\ZC^k(R^n)$ such that $\mathcal{Z}({g})\cap \zcl{X}=\zcl{X}\setminus X$. Similarly, $\nu_g$ is an isomorphism, so that $\varphi\coloneq\nu_g^{-1}\circ\nu_f$ is an isomorphism of $\ZC^k(R^n)/I$-algebras such that $\nu_f=\nu_g\circ\varphi$. It is unique by the universal property of the localization.
		\end{proof}

		\section{Cartan's theorems A and B}\label{sec:4}
		\subsection{$(X,\ZC^k_X)$ as a ringed space}
Let $X\subset R^n$ be a $\ZC$-locally closed set and consider the topology on $X$ induced by the $\ZC$-topology of $R^n$. Denote by $\text{Top}_{\ZC}(X)$ the set of open subsets of $X$ with respect to the above topology.
		
		\begin{lemma}\label{lem:sheaf}
			Let $k\in\mathbb{N}\cup\{\infty\}$. Let $X\subset R^n$ be a $\ZC$-locally closed set. Define the map
			\[
			\ZC_X^k:\text{Top}_{\ZC}(X)\to \text{Rings},\quad U\mapsto \ZC^k(U).
			\]
			The map $\ZC_X^k$ is a sheaf of rings. The couple $(X,\ZC_X^k)$ is then a ringed space.
		\end{lemma}
		
		\begin{proof}
			Consider $U,V\subset X$ two $\ZC$-open subsets such that $V\subset U$. Therefore, $U$ and $V$ are $\ZC$-locally closed subsets of $R^n$ and whenever $f\in\ZC^k(U)$, then $f|_V\in\ZC^k(V)$ by Definition \ref{def:ZC^k(X)}. Evidently, $f|_U=f$ and whenever $W\subset V$ is a $\ZC$-open subset of $X$, then $f|_W=(f|_V)|_W$. The map $\ZC_X^k$ is then a presheaf.
			
			Let $U\subset X$ be a $\ZC$-open set. Let $\{U_i\}_{i\in I}$ be a $\ZC$-open covering of $U$ in $X$ with $U_i\subset U$ for every $i\in I$. Evidently, if $f,g\in\ZC^k(U)$ satisfy $f|_{U_i}=g|_{U_i}$ for every $i\in I$, then $f=g$. Suppose then that we have functions $f_i\in\ZC^k(U_i)$ for every $i\in I$ such that $f_i|_{U_i\cap U_j}=f_j|_{U_i\cap U_j}$ for every $i,j\in I$. The function $f:U\to R$ defined as $f(x)=f_i(x)$ whenever $x\in U_i$ is well defined and it satisfies Definition \ref{def:ZC^k(X)}. Indeed, let $x\in U$, then there is $i\in I$ such that $x\in U_i$ and hence $f|_{U_i}=f_i$. As $f_i\in\ZC^k(U_i)$ there are a $\ZC$-open neighborhood $V_i$ of $x$ in $R^n$ and $F_i\in\ZC^k(V_i)$ such that $F_i|_{U_i\cap V_i}=f_i|_{U_i\cap V_i}$. As $U_i$ is $\ZC$-open in $X$ and $X$ is $\ZC$-open in $\zcl{X}$, then $\zcl{X}\setminus U_i$ is $\ZC$-closed in $R^n$ . Therefore, up to substituting $V_i$ with $V_i\setminus (\zcl{X}\setminus U_i)$, we may suppose that $U\cap V_i=U_i\cap V_i$ and then $F_i|_{U\cap V_i}=F_i|_{U_i\cap V_i}=f_i|_{U_i\cap V_i}=f|_{U_i\cap V_i}$. This shows that $f\in\ZC^k(U)$ concluding the proof that $\ZC_X^k$ is a sheaf.
		\end{proof}
		
		Moreover, we have an alternative description of the sheaf $\ZC^k_{X}$ as a direct consequence of Lemma \ref{lem:sheaf} for $X=R^n$. Let $\iota_X:X\hookrightarrow R^n$ be the inclusion map. Define the \emph{inverse image} sheaf $\iota_X^{-1}\ZC^k_{R^n}:\text{Top}_{\ZC}(X)\to \text{Rings}$ as the sheaf associated to the presheaf
		\[
		\mathcal{G}:\text{Top}_{\ZC}(X)\to \text{Rings},\quad U\mapsto \colim_{U\subset V}\ZC^k(V).
		\]
		By definition of the sheaf associated to the presheaf $\mathcal{G}$ (see \cite[\S2,\,Prop.-Def.\,1.2]{Har77}) and by Definition \ref{def:ZC^k(X)}, we directly deduce the following fact.
		
		\begin{cor}\label{cor:sheaf}
			Let $k\in\mathbb{N}\cup\{\infty\}$. Denote by $\iota_X:X\hookrightarrow R^n$ the inclusion map. Then,
			\[
			\ZC^k_X= \iota_X^{-1}\ZC^k_{R^n}/\iota_X^{-1} \I_{R^n,X}=\iota_X^{-1} \left(\ZC^k_{R^n}/\I_{R^n,X}\right),
			\]
			where $\I_{R^n,X}$ denotes the sheaf of ideals $U\mapsto \I_{\ZC^k(U)}(X\cap U)$ on $R^n$.
		\end{cor}

		\subsection{The spectrum $\operatorname{Spec}(\ZC^k(X))$ and its correspondence with $(X,\ZC^k_X)$} Let $k\in\mathbb{N}$. Denote by $\operatorname{Spec}(\ZC^k(R^n))\coloneq \{\mathfrak{p}\subset\ZC^k(R^n)\,:\, \mathfrak{p}\text{ is a prime ideal of $\ZC^k(R^n)$}\}$ the spectrum of the ring $\ZC^k(R^n)$. We denote as usual the \emph{zero locus} of $f\in\ZC^k(R^n)$ in $\operatorname{Spec}(\ZC^k(R^n))$ as
		\[
		\mathcal{V}(f)\coloneq \left\{\mathfrak{p}\in\operatorname{Spec}(\ZC^k(R^n))\,:\,  f\in\mathfrak{p}\right\}
		\]
		and, more generally, the \emph{common zero locus} of a subset $S\subset \ZC^k(R^n)$ in $\operatorname{Spec}(\ZC^k(R^n))$ as
		\[
		\mathcal{V}(S)\coloneq \left\{\mathfrak{p}\in\operatorname{Spec}(\ZC^k(R^n))\,:\,  S\subset\mathfrak{p}\right\}.
		\]
		The sets $\mathcal{V}(S)$ for $S\subset \ZC^k(R^n)$ are the closed sets of the Zariski topology of $\operatorname{Spec}(\ZC^k(R^n))$. We denote each open set for the Zariski topology of $\operatorname{Spec}(\ZC^k(R^n))$ as
		\[
		\mathcal{U}(S)\coloneq \operatorname{Spec}\left(\ZC^k(R^n)\right)\setminus \mathcal{V}(S)=\left\{\mathfrak{p}\in\operatorname{Spec}(\ZC^k(R^n))\,:\,  S\nsubset\mathfrak{p}\right\}.
		\]
		
		Consider the inclusion map $\jmath_{R^n}:R^n\to \operatorname{Spec}(\ZC^k(R^n))$ defined as $\jmath_{R^n}(a)\coloneq \mathfrak{m}_{a}$. Observe that the map $\jmath_{R^n}:R^n\to \operatorname{Spec}\left(\ZC^k(R^n)\right)$ is injective and it is continuous with respect to the $\ZC$-topology of $R^n$ and the Zariski topology of $\operatorname{Spec}\left(\ZC^k(R^n)\right)$. Moreover, by Corollary \ref{cor:max-ideals}, we have that $\jmath_{R^n}(R^n)=\operatorname{Spec}_m\left(\ZC^k(R^n)\right)$, the set of maximal ideals of $\ZC^k(R^n)$.\\

		Let $Y\subset R^n$ be a $\ZC$-closed subset. Denote by
		\[
		\widetilde{Y}\coloneq \mathcal{V}(\jmath_{R^n}(Y))
		\]
		the Zariski closure of $\jmath_{R^n}(Y)$ in $\operatorname{Spec}(\ZC^k(R^n))$. Observe that $\mathcal{Z}(\I_{\ZC^k(R^n)}(Y))=Y$ for every $\ZC$-closed subset $Y$ of $R^n$; hence, whenever $Z,Y\subset R^n$ are $\ZC$-closed subsets, $Z\subset Y$ if and only if $\widetilde{Z}\subset\widetilde{Y}$.
		
		As a consequence of the Nullstellensatz for $\ZC^k(R^n)$ (see Theorem \ref{thm:Nullstellensatz}) and of Theorem \ref{thm:radprinc}, we have the following characterization of Zariski closed subsets of $\operatorname{Spec}(\ZC^k(R^n))$.
		
		\begin{thm}\label{thm:correspondence}
			Let $k\in \mathbb{N}$. Every Zariski closed subset of $\operatorname{Spec}(\ZC^k(R^n))$ is of the form $\mathcal{V}(g)$ for some $g\in \ZC^k(R^n)$.
			\vspace{0.5em}
			
			Moreover, $\mathcal{V}(g_1)=\mathcal{V}(g_2)$ in $\operatorname{Spec}(\ZC^k(R^n))$ if and only if $\mathcal{Z}(g_1)=\mathcal{Z}(g_2)$ in $R^n$. In particular, there is a correspondence between $\ZC$-closed subsets $Y=\mathcal{Z}(g)$ of $R^n$ and the unique Zariski closed set $\widetilde{Y}\coloneq \mathcal{V}(\jmath_{R^n}(Y))=\mathcal{V}(g)$ of $\operatorname{Spec}(\ZC^k(R^n))$ such that $\jmath_{R^n}(Y)=\widetilde{Y}\cap\operatorname{Spec}_m(\ZC^k(R^n))$.
		\end{thm}
		\begin{proof}
			Let $I\subset\ZC^k(R^n)$ be an ideal. Recall that $\ZC^k(R^n)$ is radically principal (see Theorem \ref{thm:radprinc}) and that, by Corollary \ref{cor:radgeneration&zeros}, $\sqrt{I}=\sqrt{(g)}$ for $g\in I$ if and only if $\mathcal{Z}(g)=\mathcal{Z}(I)$. Therefore, $\mathcal{V}(I)=\mathcal{V}\left(\sqrt{I}\right)=\mathcal{V}\left(\sqrt{(g)}\right)=\mathcal{V}(g)$. Moreover,
			$\mathcal{V}(g_1)=\mathcal{V}(g_2)$ for $g_1,g_2\in\ZC^k(R^n)$ if and only if $\sqrt{(g_1)}=\sqrt{(g_2)}$, that is, if and only if $\mathcal{Z}(g_1)=\mathcal{Z}(g_2)$ again by Corollary \ref{cor:radgeneration&zeros}. This shows that the map $Y\mapsto \widetilde{Y}\coloneq \mathcal{V}(\jmath_{R^n}(Y))$, where $Y$ ranges on the set of $\ZC$-closed subsets of $R^n$, is well defined and it is bijective onto the Zariski closed subsets of $\operatorname{Spec}\left(\ZC^k(R^n)\right)$. In addition, if $g\in\ZC^k(R^n)$ satisfies $Y=\mathcal{Z}(g)$, then
			\begin{align*}
				\widetilde{Y}\cap\operatorname{Spec}_m\left(\ZC^k(R^n)\right)
				&=\left\{\mathfrak{m}_{a}\in\operatorname{Spec}_m\left(\ZC^k(R^n)\right)\,:\, a\in\mathcal{Z}(g)\right\}=\jmath_{R^n}(\mathcal{Z}(g))=\jmath_{R^n}{(Y)},
			\end{align*}
			as desired.
		\end{proof}
		
		\begin{rem}
			Observe that Theorem \ref{thm:correspondence} does not hold when replacing the ring $\ZC^k(R^n)$ with the usual ring of regular functions $\mathcal{R}(R^n)$ on $R^n$. Indeed, set $n=2$, $Y\coloneq \{(0,0)\}$ and denote by $\jmath:R^2\hookrightarrow\operatorname{Spec}(\mathcal{R}(R^2))$ the inclusion $\jmath(a)\coloneq (x_1-a_1,x_2-a_2)$ for every $a=(a_1,a_2)\in R^2$. Observe that both $\mathcal{V}\left(\{x,y\}\right)$ and $\mathcal{V}\left(x^2+y^2\right)$ in $\operatorname{Spec}\left(\mathcal{R}(R^2)\right)$ satisfy $\jmath^{-1}(\mathcal{V}(\{x,y\}))=Y=\jmath^{-1}\left(\mathcal{V}\left(x^2+y^2\right)\right)$ but $\mathcal{V}(\{x,y\})$ is not principal and $\mathcal{V}(\{x,y\})\subsetneq\mathcal{V}\left(x^2+y^2\right)$ as $\left(x^2+y^2\right)$ is a prime ideal in $\mathcal{R}(R^2)$ strictly contained in $(x,y)$. By contrast, Corollary \ref{cor:radgeneration&zeros} ensures that $\sqrt{(x,y)\cdot\ZC^k(R^2)}=\sqrt{\left(x^2+y^2\right)\cdot\ZC^k(R^2)}$ and hence $\mathcal{V}(\{x,y\})=\mathcal{V}\left(x^2+y^2\right)$ in $\operatorname{Spec}\left(\ZC^k(R^2)\right)$.\\
		\end{rem}

		Let $k\in\mathbb{N}$. Denote by $\widetilde{R^n}\coloneq \mathcal{V}(\jmath_{R^n}(R^n))=\operatorname{Spec}\left(\ZC^k(R^n)\right)$. Then, by Theorem \ref{thm:correspondence}, the usual structure sheaf $\mathcal{O}_{\widetilde{R^n}}$ on $\widetilde{R^n}$ has a very simple description: let $\mathcal{U}(S)$ be a Zariski open subset of $\widetilde{R^n}$, then there is $g\in\ZC^k(R^n)$ such that $\mathcal{U}(g)=\mathcal{U}(S)$ and therefore
		\[
		\mathcal{O}_{\widetilde{R^n}}(\mathcal{U}(S))=\mathcal{O}_{\widetilde{R^n}}(\mathcal{U}(g))\coloneq \ZC^k(R^n)_g.
		\]
		
		Consider the inclusion map $\jmath_{R^n}:R^n\hookrightarrow \operatorname{Spec}(\ZC^k(R^n))$. By Theorem \ref{thm:correspondence}, for every $\ZC$-open set $U$ of $R^n$ there is a unique Zariski open subset $\widetilde{U}$ of $\operatorname{Spec}(\ZC^k(R^n))$ such that $\jmath_{R^n}(U)=\widetilde{U}\cap\operatorname{Spec}_m(\ZC^k(R^n))$. Therefore, the inverse image sheaf $\jmath_{R^n}^{-1} \mathcal{O}_{\widetilde{R^n}}$ satisfies 
		\[
		\jmath_{R^n}^{-1}\mathcal{O}_{\widetilde{R^n}}(U)=\mathcal{O}_{\widetilde{R^n}}\left(\widetilde{U}\right)\coloneq \ZC^k(R^n)_g
		\]
		for every $\ZC$-open subset $U$ of $R^n$ with $g\in\ZC^k(R^n)$ satisfying $\mathcal{Z}(g)=R^n\setminus U$.\\
		
		Let $X\subset R^n$ be a $\ZC$-locally closed set. Denote by $\I_{\widetilde{R^n},X}$ the sheaf of ideals $\mathcal{U}(g)\mapsto \I_{\ZC^k(R^n)}(X)\cdot \ZC^k(R^n)_g$ on $\operatorname{Spec}\left(\ZC^k(R^n)\right)$. Observe that, by Proposition \ref{prop:localization}, we get that:
		\begin{equation}\label{eq:sheaves}
			\iota_{X}^{-1}\left(\jmath_{R^n}^{-1}\left(\mathcal{O}_{\widetilde{R^n}}/\I_{\widetilde{R^n},X}\right)\right)\cong\iota_{X}^{-1}\left(\ZC^k_{R^n}/\I_{R^n,X}\right)=\ZC^k_X.
		\end{equation}
		
		Let us look for another characterization of the sheaf $\ZC^k_X$, which we will obtain also by using (\ref{eq:sheaves}).

	Let $X\subset R^n$ be a $\ZC$-locally closed set. Let $f\in\ZC^k(R^n)$ such that $\mathcal{Z}(f)\cap\zcl{X}=\zcl{X}\setminus X$. Denote by
	\[
	\widetilde{X}\coloneq\operatorname{Spec}\left(\left(\ZC^k(R^n)/\I_{\ZC^k(R^n)}(X)\right)_{\bar{f}}\right).
	\]
	Consider the inclusion $\iota_{\widetilde{X}}:\widetilde{X}\hookrightarrow\operatorname{Spec}\left(\ZC^k(R^n)\right)$ defined as $\iota_{\widetilde{X}}(\mathfrak{p})=\pi^{-1}\left(\mathfrak{p}\cap\ZC^k(R^n)/\I_{\ZC^k(R^n)}(X)\right)$ for every $\mathfrak{p}\in\widetilde{X}$, where $\pi:\ZC^k(R^n)\to\ZC^k(R^n)/\I_{\ZC^k(R^n)}(X)$ denotes the quotient morphism. Observe that $\iota_{\widetilde{X}}(\widetilde{X})=\mathcal{V}\left(\I_{\ZC^k(R^n)}(X)\right)\cap \mathcal{U}(f)$. Endow $\widetilde{X}$ with its Zariski topology, which coincides with the induced one by the Zariski topology of $\operatorname{Spec}\left(\ZC^k(R^n)\right)$.

Denote by $\jmath_X:X\hookrightarrow \widetilde{X}$ the inclusion defined by
\[
\jmath_X(a)\coloneq \left(\mathfrak{m}_a/\mathcal{I}_{\ZC^k(R^n)}(X)\right)\cdot \left(\ZC^k(R^n)/\mathcal{I}_{\ZC^k(R^n)}(X)\right)_{\bar{f}}.
\]
Observe that, by standard properties of quotient and localization of rings, we have that
\[
\jmath_X(X)=\operatorname{Spec}_m\left(\left(\ZC^k(R^n)/\mathcal{I}_{\ZC^k(R^n)}(X)\right)_{\bar{f}}\right).
\]
Denote by $\mathcal{O}_{\widetilde{X}}$ the structure sheaf of $\widetilde{X}$. Observe that:
\begin{equation}\label{eq:spec-loc}
	\mathcal{O}_{\widetilde{X}}=\iota_{\widetilde{X}}^{-1} 	\left(\mathcal{O}_{\widetilde{R^n}}/\mathcal{I}_{\widetilde{R^n},X}\right),
\end{equation}
where $\mathcal{I}_{\widetilde{R^n},X}$ denotes the sheaf of ideals $\mathcal{U}(g)\mapsto \mathcal{I}_{\ZC^k(R^n)}(X)\cdot \ZC^k(R^n)_g$ on $\operatorname{Spec}\left(\ZC^k(R^n)\right)$.\\

\begin{rem}\label{rem:spec}
	Denote by
	\[
	\nu_f:\left(\ZC^k(R^n)/\I_{\ZC^k(R^n)}(X)\right)_{\bar{f}}\to\ZC^k(X)
	\]
	the (canonical) isomorphism of Corollary \ref{cor:extension}. It induces an isomorphism of affine schemes
	\[
	\widetilde{\nu_f}:\left(\operatorname{Spec}\left(\ZC^k(X)\right),\mathcal{O}_{\operatorname{Spec}\left(\ZC^k(X)\right)}\right)\to\left(\widetilde{X},\mathcal{O}_{\widetilde{X}}\right),\quad\mathfrak{q}\mapsto\nu_f^{-1}(\mathfrak{q}).
	\]
\end{rem}
		
		Therefore, let us simplify the notation by identifying $\operatorname{Spec}\left(\left(\ZC^k(R^n)/\I_{\ZC^k(R^n)}(X)\right)_{\bar{f}}\right)$ with $\operatorname{Spec}\left(\ZC^k(X)\right)$ as follows:
		
		\begin{notation}\label{not:2}
			We settle the following substitutions: ``$\widetilde{X}$"$\coloneq \operatorname{Spec}\left(\ZC^k(X)\right)$, ``$\mathcal{O}_{\widetilde{X}}$"$\coloneq \mathcal{O}_{\operatorname{Spec}\left(\ZC^k(X)\right)}$, ``$\jmath_X$"$\coloneq \widetilde{\nu_f}^{-1}\circ\jmath_X$ and ``$\iota_{\widetilde{X}}$"$\coloneq \iota_{\widetilde{X}}\circ\widetilde{\nu_f}$.
		\end{notation}
		
		In light of Notation \ref{not:2}, the following diagram of inclusions commutes:
		\begin{equation}\label{diagram}
			\begin{tikzcd}
				R^n \arrow[r,hook, "\jmath_{R^n}"]  & \operatorname{Spec}\left(\ZC^k(R^n)\right) \\
				X \arrow[u,hook, "\iota_X"] \arrow[r,hook, "\jmath_X"] & \operatorname{Spec}\left(\ZC^k(X)\right) \arrow[u,hook,"\iota_{\widetilde{X}}"]
			\end{tikzcd}
		\end{equation}
		
		Let $Y\subset X$ be a $\ZC$-closed subset. Denote by
		\[
		\widetilde{Y}\coloneq \mathcal{V}(\jmath_{X}(Y))
		\]
		the Zariski closure of $\jmath_{X}(Y)$ in $\operatorname{Spec}\left(\ZC^k(X)\right)$. Observe that, by Proposition \ref{prop:zerolocus}, if $Z,Y\subset X$ are $\ZC$-closed subsets, then $Z\subset Y$ if and only if $\widetilde{Z}\subset\widetilde{Y}$. Moreover, $\widetilde{X}\coloneq\mathcal{V}(\jmath_{X}(X))=\mathcal{V}(0)=\operatorname{Spec}\left(\ZC^k(X)\right)$.\\
		
		Let $k\in\mathbb N$. As a consequence of Theorem \ref{thm:correspondence} and the usual properties of quotients and localization of ideals, we have the following characterization of Zariski closed subsets of $\operatorname{Spec}\left(\ZC^k(X)\right)$.
		
		\begin{cor}\label{cor:correspondence}
			Let $k\in \mathbb{N}$. Let $X\subset R^n$ be a $\ZC$-locally closed set. Every Zariski closed subset of $\operatorname{Spec}\left(\ZC^k(X)\right)$ is of the form $\mathcal{V}(g)$ for some $g\in \ZC^k(X)$.
			\vspace{0.5em}
			
			Moreover, $\mathcal{V}({g}_1)=\mathcal{V}({g}_2)\subset\operatorname{Spec}\left(\ZC^k(X)\right)$ if and only if $\mathcal{Z}(g_1)=\mathcal{Z}(g_2)\subset X$. In particular, there is a correspondence between $\ZC$-closed subsets $Y=\mathcal{Z}(g)$ of $X$ and the unique Zariski closed set $\widetilde{Y}\coloneq \mathcal{V}(\jmath_{X}(Y))=\mathcal{V}({g})$ of $\operatorname{Spec}\left(\ZC^k(X)\right)$ such that $\jmath_{X}(Y)=\widetilde{Y}\cap\operatorname{Spec}_m\left(\ZC^k(X)\right)$.
		\end{cor}
		
		By commutativity of the diagram (\ref{diagram}) and by equations (\ref{eq:sheaves})\&(\ref{eq:spec-loc}), we get that:
		\[
		\ZC^k_X\cong(\iota_{\widetilde{X}}\circ \jmath_{X})^{-1}\left(\mathcal{O}_{\widetilde{R^n}}/\I_{\widetilde{R^n},X}\right)=\jmath^{-1}_X \mathcal{O}_{\widetilde{X}}.
		\]
		
		Let us summarize the previous results as follows.
		
		\begin{cor}\label{cor:inverse-sheaf}
			Let $k\in \mathbb{N}$. Let $X\subset R^n$ be a $\ZC$-locally closed set. Let $\jmath_X:X\to \operatorname{Spec}\left(\ZC^k(X)\right)$ be the inclusion map defined as above. Then:
			\[
			\ZC^k_X\cong \jmath_X^{-1} \mathcal{O}_{\widetilde{X}},
			\]
			where $\widetilde{X}\coloneq\operatorname{Spec}\left(\ZC^k(X)\right)$ and $\mathcal{O}_{\widetilde{X}}$ denotes its structure sheaf.
		\end{cor}

		\subsection{Sheaves of modules on $X$ and $\operatorname{Spec}(\ZC^k(X))$: Cartan's theorems A and B}\label{sub:Cartan}
		
		Let $k\in\mathbb{N}$. Let $X\subset R^n$ be a $\ZC$-locally closed set and denote by $\widetilde{X}=\operatorname{Spec}\left(\ZC^k(X)\right)$. Let us fix some standard notation for ringed spaces in the special case of $\left(X,\ZC^k_X\right)$ and $\left(\widetilde{X},\mathcal{O}_{\widetilde{X}}\right)$.\\
		
		Let $M$ be a $\ZC^k(X)$-module. Let us recall the construction of the associated $\ZC^k_X$-module $\widetilde{M}_{X}$ on $X$. Let $U\subset X$ be a $\ZC$-open subset of $X$ and let $g\in\ZC^k(X)$ such that $\mathcal{Z}(g)=X\setminus U$ (see Proposition \ref{prop:zerolocus} and Corollary \ref{cor:extension}). Then, define:
		\[
		\widetilde{M}_X(U)\coloneq M_{g}.
		\]
        Note that if $g'\in\ZC^k(X)$ satisfies $\mathcal{Z}(g')=X\setminus U$, then, as discussed in Corollary \ref{cor:correspondence}, we deduce that $\sqrt{(g)}=\sqrt{(g')}$ in $\ZC^k(X)$. Therefore, $g'^N\in(g)$ and $g^N\in(g')$ for some $N\in\mathbb{N}\setminus\{0\}$, and hence  $M_g\simeq M_{g'}$.
		
		Let $M$ be a $\ZC^k(X)$-module. Let us recall the construction of the associated $\mathcal{O}_{\widetilde{X}}$-module $\widetilde{M}_{\widetilde{X}}$ on $\widetilde{X}$. By Corollary \ref{cor:correspondence} and Remark \ref{rem:spec}, every Zariski open subset of $\widetilde{X}$ is of the form $\mathcal{U}(g)\coloneq \widetilde{X}\setminus \mathcal{V}(g)$, for some $g\in\ZC^k(X)$. Then, define:
		\[
		\widetilde{M}_{\widetilde{X}}\left(\mathcal{U}({g})\right)\coloneq M_g.
		\]
		
		Observe that by Corollary \ref{cor:correspondence} and the very definitions, we have that
		\[
		\widetilde{M}_{X}=\iota_{X}^{-1} \widetilde{M}_{\widetilde{X}}\quad\text{and}\quad\widetilde{M}_{\widetilde{X}}=(\iota_{X})_* \widetilde{M}_{X},
		\]
		and, more generally, Corollaries \ref{cor:correspondence}\&\ref{cor:inverse-sheaf} ensure that the maps $\iota_{X}^{-1}$ and $(\iota_{X})_*$ are bijections between the sets of $\ZC^k_X$-modules on $X$ and the set of  $\mathcal{O}_{\widetilde{X}}$-modules on $\widetilde{X}$. Moreover, these operators are inverse of each other, that is, $\iota_{X}^{-1}\circ(\iota_{X})_*(\mathcal{F})=\mathcal{F}$ for every $\ZC^k_X$-module $\mathcal{F}$ on $X$ and $(\iota_{X})_*\circ\iota_{X}^{-1}(\mathcal{G})=\mathcal{G}$ for every $\mathcal{O}_{\widetilde{X}}$-module $\mathcal{G}$ on $\widetilde{X}$.\\
		
		Let us recall a few general notions for completeness. Let $(Z,\mathcal{O}_Z)$ be a ringed space. A $\mathcal{O}_Z$-module $\mathcal{F}$ on $Z$ is said to be \emph{quasi-coherent} if there is an open cover $\{U_\alpha\}$ of $Z$ such that for every $\alpha$ there is a $\mathcal{O}(U_\alpha)$-module $M_\alpha$ such that $\mathcal{F}|_{U_\alpha}\cong \widetilde{M_\alpha}$, where $ \widetilde{M_\alpha}$ denotes the $\mathcal{O}(U_\alpha)$-module on $U_\alpha$ associated to $M_\alpha$. A $\mathcal{O}_Z$-module $\mathcal{F}$ on $Z$ is said to be \emph{generated by its global sections} if for every $z\in Z$, for every open subset $U$ of $Z$ and for every $\sigma\in\mathcal{F}(U)$ there are:\\
		
		\begin{enumerate}[label=(\roman*),nosep]
			\item an open neighborhood $V$ of $z$ in $U$;
			\item $\ell\in \mathbb{N}$, $\sigma_1,\dots,\sigma_\ell\in\mathcal{F}(Z)$ and $a_1,\dots,a_\ell\in\mathcal{O}_Z(V)$ such that $\sigma|_V=\sum_{i=1}^{\ell} a_i \sigma_i|_{V}$.\\
		\end{enumerate}
		
		Equivalently, the last property can be restated in terms of stalks: for every $z\in Z$ the image of $\mathcal{F}(Z)$ in $\mathcal{F}_z$ generates $\mathcal{F}_z$ as an $\mathcal{O}_{Z,z}$-module. \\
		
		\begin{rem}\label{rem:modules}
			Let $\mathcal{F}$ be a sheaf of $\ZC^k_X$-modules on $X$ and consider the sheaf $(\iota_{X})_*\mathcal{F}$ of $\mathcal{O}_{\widetilde{X}}$-modules on $\widetilde{X}$. By Corollary \ref{cor:correspondence} we have that:
			\begin{enumerate}[label=(\roman*),nosep]
				\item\label{rem:1} $\mathcal{F}$ is quasi-coherent if and only if $(\iota_{X})_*\mathcal{F}$ is so.
				\item\label{rem:2} $\mathcal{F}$ is generated by its global sections if and only if $(\iota_{X})_*\mathcal{F}$ is so.
				\item\label{rem:3} $H^p(X,\mathcal{F})=0$ for every $p>0$ if and only if $H^p\left(\widetilde{X},(\iota_{X})_*\mathcal{F}\right)=0$ for every $p>0$.\\
			\end{enumerate} 
		\end{rem}
		
		We have now all the information to deduce the following version of Cartan's Theorem A and B (see \cite[p.\,51]{Car53}) in the algebraic formulation of Serre (see \cite[\S45,Thm.\,2;\,\S46,Thm.\,3]{Ser55}). As $\ZC^k(X)$ is not Noetherian, an application of \cite[Thm.\,1.3.1]{Gro61} together with Remark \ref{rem:modules}, gives the following result.
		
		\begin{thm}[Cartan's Theorems A and B]\label{thm:Cartan-A&B}
			Let $k\in\mathbb{N}$. Let $X\subset R^n$ be a $\ZC$-locally closed set and let $\mathcal{F}$ be a quasi-coherent $\ZC^k_X$-module on $X$. Then:
			
			\begin{enumerate}[label=(\Alph*),nosep]
				\item\label{A} $\mathcal{F}$ is generated by its global sections.
				\item\label{B} $H^p(X,\mathcal{F})=0$ for every $p>0$.
			\end{enumerate} 
		\end{thm}
		
		As a direct consequence of Theorem \ref{thm:Cartan-A&B}\ref{B} we get an alternative proof of the extension of $\ZC^k$-functions on a $\ZC$-locally closed subset $X$ of $R^n$ to any $\ZC$-open neighborhood $U$ of $X$ in $R^n$ such that $X=\zcl{X}\cap U$.
		
		\begin{proof}[Proof of Theorem \ref{thm:extension}] Let $X\subset R^n$ be a $\ZC$-locally closed set and let $U\subset R^n$ be a $\ZC$-open set such that $X=U\cap\zcl{X}$. Denote by $\iota_{U}:U\hookrightarrow R^n$ and $\iota_{X}:X\hookrightarrow R^n$ the inclusion maps. Consider the following exact sequence of coherent sheaves on $U$:
		\[
			\begin{tikzcd}
				0\arrow[r] & \iota_U^{-1} \I_{R^n,X} \arrow[r]  & \iota_U^{-1} \ZC^k_{R^n} \arrow[r] & \iota_U^{-1} \ZC^k_{R^n}/\iota_U^{-1}\I_{R^n,X} \arrow[r] & 0.
			\end{tikzcd}
		\]
		As $\iota_U^{-1} \ZC^k_{R^n}=\ZC^k_U$, Theorem \ref{thm:Cartan-A&B}\ref{B} ensures that the induced sequence
		\[
		\begin{tikzcd}
			0\arrow[r] &\I_{\ZC^k(U)}(X) \arrow[r]  & \ZC^k(U) \arrow[r] & \iota_U^{-1} \left(\ZC^k_{R^n}/\I_{R^n,X}\right)(U) \arrow[r] & 0
		\end{tikzcd}
		\]
		is exact. As the sheaf $\ZC^k_{R^n}/\I_{R^n,X}$ is supported over $\zcl{X}$, we get that $\iota_U^{-1} \left(\ZC^k_{R^n}/\I_{R^n,X}\right)(U)=\iota_X^{-1} \left(\ZC^k_{R^n}/\I_{R^n,X}\right)(X)=\ZC^k(X)$ by Corollary \ref{cor:sheaf}. This shows that the restriction morphism $\rho:\ZC^k(U)\to \ZC^k(X)$, $F\mapsto F|_{U\cap X}$ is surjective, as desired.
		\end{proof}
		
		\printbibliography[heading=bibintoc]
		
	\end{document}